\documentclass[10pt,a4paper,sumlimits,namelimits]{article}

\usepackage{amsthm}
\usepackage{amsmath}
\usepackage{amsfonts}
\usepackage{amssymb}
\usepackage{mathrsfs}
\usepackage{empheq}
\usepackage{stmaryrd,soul}
\usepackage{enumerate}
\usepackage[applemac]{inputenc}
\usepackage[colorlinks=true,urlcolor=blue]{hyperref}

\newcommand{\N}{\mathbb{N}}

\newcommand{\R}{\mathbb{R}}
\newcommand{\C}{\mathbb{C}}

\newcommand{\Dr}{\mathscr{D}}
\newcommand{\Lr}{\mathscr{L}}
\newcommand{\vphi}{\varphi}
\newcommand{\eps}{\varepsilon}

\newcommand{\dsp}{\displaystyle}
\newcommand{\ovl}{\overline}
\newcommand{\udl}{\underline}
\newcommand{\vlim}{\lim\limits}

\newcommand{\vinf}{\inf\limits}
\newcommand{\vint}{\int\limits}

\newcommand{\inj}{\hookrightarrow}
\newcommand{\tends}{\longrightarrow}
\newcommand{\weak}{\rightharpoonup}

\newcommand{\loc}{\mathrm{loc}}

\renewcommand{\b}{\mathrm{b}}
\newcommand{\co}{\mathrm{c}}
\renewcommand{\d}{\mathrm{d}}
\newcommand{\vu}{\mathrm{u}}

\newcommand{\GN}{\mathrm{GN}}
\newcommand{\vi}{\mathrm{i}}
\newcommand{\w}{{\textsl w}}
\renewcommand{\le}{\leqslant}
\renewcommand{\ge}{\geqslant}
\renewcommand{\Re}{\mathrm{Re}}
\renewcommand{\Im}{\mathrm{Im}}
\newcommand{\bs}{\boldsymbol}
\newcommand{\p}{\prime}

\DeclareMathOperator{\supp}{supp}

\DeclareMathOperator{\argmax}{argmax}

\numberwithin{equation}{section}

\newtheorem{thm}{Theorem}[section]
\newtheorem{prop}[thm]{Proposition}

\newtheorem{lem}[thm]{Lemma}

\theoremstyle{definition}
\newtheorem{rmk}[thm]{Remark}
\newtheorem{defi}[thm]{Definition}

\newtheorem{assum}[thm]{Assumption}

\newenvironment{proof*}{\noindent{\bf Proof.}}{\qed}
\newenvironment{vproof}[1]{\noindent{\bf Proof #1}}{\qed}

\title{\huge \sc Finite Time Extinction for the Strongly Damped Nonlinear Schr\"{o}dinger Equation in Bounded Domains}
\author{\sc Pascal Bégout$^*$ and Jes\'us Ildefonso D{\'{\i}}az$^\dagger$}
\date{}

\begin{document}

\maketitle

\begin{gather*}
\begin{array}{cc}
            ^*\text{Institut de Mathématiques de Toulouse \& TSE }	&	^\dagger\text{Instituto de Matem\'atica
																			Interdisciplinar}	\\
                                         \text{Université Toulouse I Capitole }	&	\text{ Departamento de An\'alisis y Matem\'atica
                                         															Aplicada}		\\
                                                  \text{Manufacture des Tabacs }	&	\text{ Universidad Complutense de Madrid}	\\
                                                         \text{21, Allée de Brienne }	&	\text{ Plaza de las Ciencias, 3}				\\
                                 \text{31015 Toulouse Cedex 6, FRANCE }	&	\text{ 28040 Madrid, SPAIN}
\bigskip \\
\text{
{\footnotesize E-mail\:: }\htmladdnormallink{{\footnotesize\udl{\tt{Pascal.Begout@math.cnrs.fr}}}}{mailto:Pascal.Begout@math.cnrs.fr}}
&
\text{
{\footnotesize E-mail\:: }\htmladdnormallink{{\footnotesize\udl{\tt{jidiaz@ucm.es}}}}{mailto:jidiaz@ucm.es}
}
\end{array}
\end{gather*}

\begin{abstract}
We prove the \textit{finite time extinction property} $(u(t)\equiv 0$ on $\Omega$ for any $t\ge T_\star,$ for some $T_\star>0)$ for solutions of the nonlinear Schr\"{o}dinger problem $\vi u_t+\Delta u+a|u|^{-(1-m)}u=f(t,x),$ on a bounded
domain $\Omega$ of $\R^N,$ $N\le 3,$ $a\in\C$ with $\Im(a)>0$ (the damping case) and under the crucial assumptions $0<m<1$ and the dominating condition $2\sqrt m\,\Im(a)\ge(1-m)|\Re(a)|.$ We use an energy method as well as several a priori estimates to prove the main conclusion. The presence of the non-Lipschitz nonlinear term in the equation introduces a lack of regularity of the solution requiring a study of the existence and uniqueness of solutions satisfying the equation in some different senses according to the regularity assumed on the data. 
\end{abstract}

{\let\thefootnote\relax}
{\let\thefootnote\relax}

\tableofcontents

\baselineskip .6cm

\section{Introduction}
\label{introduction}

This paper deals with the \textit{finite time extinction property} of solutions of the nonlinear Schr\"{o}dinger problem
\begin{gather}
\begin{cases}
\label{nlsi}
	\vi\dfrac{\partial u}{\partial t}+\Delta u+a|u|^{-(1-m)}u=f(t,x),	&	\text{in } (0,\infty)\times\Omega,		\medskip \\
	u(t)_{|\Gamma}=0,									&	\text{on } (0,\infty)\times\Gamma,	\medskip \\
	u(0)= u_0,											&	\text{in } \Omega,
\end{cases}
\end{gather}
when, roughly speaking, we assume that $N\le3,$
\begin{equation}
\label{Damping}
a\in\C \; \text{ with } \; \Im (a)>0,
\end{equation}
and 
\begin{equation}
\label{strongly}
0<m<1.
\end{equation}
We start by pointing out that this \textit{finite time extinction property} $(u(t)\equiv 0$ on $\Omega$ for any $t\ge T_\star,$ for some $T_\star>0)$ represents, clearly, the most opposite property to the famous Max Born result on the \textit{conservation of the mass} 
\begin{gather*}
\|u(t)\|_{L^2(\Omega )}=\|u_0\|_{L^2(\Omega)}, \text{ for any } t\ge0,
\end{gather*}
which arises (when $f=0)$ in the linear case (and more generally if $\Im(a)=0:$ see Proposition~\ref{propL2} below) and which allows the probabilistic understanding of the complex wave solution $u(t,x)$ in the context of the applications of the linear Schr\"{o}dinger equation in Quantum Mechanics. It is well known that the presence of a damping term \eqref{Damping} makes
the equation irreversible with respect the time.

We also recall that the Schr\"{o}dinger equation in presence of a nonlinear term in the equation (as, e.g., problem~\eqref{nlsi} when $a\in\C$ and $a\neq 0)$ arises in many other different contexts as, e.g., Nonlinear Optics, Hydrodynamics, etc., and that those other contexts, for instance in Nonlinear Optics, the variable $t$ does not represent time but the main scalar spacial variable which appears in the propagation of the waveguide direction (see e.g. Agrawal and Kivshar~\cite{ak}, Sulem and Sulem~\cite{MR2000f:35139}, Shi, Xu, Yang, Yang and Yin~\cite{MR3253582} and its many references).

As a matter of fact, the nonlinear Schr\"{o}dinger equation under condition \eqref{Damping} is referred in the literature as the \textit{damped} case and it was intensively studied since the middle of the past century under different additional conditions (but most of them for $m>1)$ (see, e.g., Nelson~\cite{MR0161189}, Pozzi~\cite{MR0229987}, Bardos and Brezis~\cite{MR0242020}, Lions~\cite{MR0259693}, Kato~\cite{MR0492961}, Brezis and Kato~\cite{MR80i:35135}, Vladimirov~\cite{MR745511}, Tsutsumi~\cite{MR1038160}, Temam and Miranville~\cite{MR2169020}, Kita and Shimomura~\cite{MR2272871}, Carles and Gallo~\cite{MR2765425}, Carles and Ozawa~\cite{MR3306821} and Hayashi, Li and Naumkin~\cite{MR3465033}, among others).

In our above formulation we assume that $a\in \mathbb{C}$ and thus a possible, non-dominant non-dissipative nonlinear term may coexists with the damping term (i.e., we allow $\Re (a)\neq 0).$ Nevertheless, our main result on the finite time extinction for $|\Omega |<\infty$ requires the dominating condition
\begin{gather*}
2\sqrt m\,\Im(a)\ge(1-m)|\Re(a)|,
\end{gather*}
as well as the assumption~\eqref{strongly} on a strong damping. 

We also recall that in most of the papers on the nonlinear equation \eqref{nlsi} it is assumed that $m=3$ (the so called cubic case). Nevertheless there are several applications in which the general case $m>0$ is of interest. For instance, it is the case of the so called \textit{non-Kerr type equations} arising in the study of optical solitons (see, e.g., \cite{ak}). For
some other physical details and many references, we refer the reader to the general presentations made in the books \cite{ak} and \cite{MR2000f:35139}. Some other references concerning the case $m\in (0,1)$ are quoted in our previous paper Bégout and D\'{\i}az~\cite{MR3193996}. We also mention that the spacial localization phenomenon (solutions with support $u(t,\:.\:)$ being a compact, when $\Omega $ is unbounded) requires a different balance between the damping and non-damping components (mainly with $\Im(a)>0$) of the nonlinear term $a|u|^{m-1}u$ (see \cite{MR2876246,MR3193996,MR3190983}).

In spite of the large amount of papers devoted to the existence and uniqueness results of nonlinear Schr\"{o}dinger equations with a damping term only very few of them allowed the consideration of a strong damping term (i.e. condition \eqref{strongly}). This is the reason why we presented here some new results on the general theory of the existence, uniqueness and
regularity of solutions of the strongly damped Schr\"odinger equation improving several previous papers in the literature (see, e.g. Carles and Gallo~\cite{MR2765425}, Lions~\cite{MR0259693}, Brezis and Cazenave~\cite{bc} and Vrabie~\cite{MR1375237}) which are needed for the study of the finite time extinction property.

Since the comparison principle does not apply to our problem, the main tool to prove the finite time extinction property is a suitable \textit{energy method} in the spirit of the collection of energy methods quoted in the monograph Antontsev, D\'{\i}az and Shmarev~\cite{MR2002i:35001}. Nevertheless, the adaptation to the nonlinear Schr\"{o}dinger equation requires some new estimates and also a sharper study of the ordinary differential inequality satisfied by the mass. We start by giving, in Section~\ref{general}, a \textit{semi-abstract} result (which is proved in Section~\ref{prooffinite}) in which the finite time extinction property is derived under a general regularity condition on the solution. The presence of the non-Lipschitz nonlinear term in the equation introduces a lack of regularity of the solution (in contrast to the case in which $m\ge1)$ and so we shall devote Section~\ref{exiuni} to present a separated study of the existence and uniqueness of solutions satisfying the equation in some different sense according to the regularity assumed on the data. To this purpose, we use mainly some monotonicity methods, jointly with suitable regularizations and passing to the limit, improving previous results in the literature. Section~\ref{finite} concerns the finite time extinction and the asymptotic behavior of the solution. The proofs of the results of Sections~\ref{finite} and \ref{exiuni} are presented in Sections~\ref{proofext} and \ref{proofexi}, respectively. An Appendix, collecting some technical auxiliary results, is also presented for the convenience of the reader.

We point out that in our formulation it may arise a non-homogeneous term (on which we assume a finite time extinction $T_0)$ and that, surprisingly enough, under some critical decay to zero of $f(t,\:.\:)$ at $t=T_0,$ we can conclude that the corresponding solution $u$ also vanishes after the same time $t=T_0$ (see Theorem~\ref{thmgenext} part~\ref{pthmgenext2}). Our energy method allows us also to get some large time decay estimates in some cases, always under the presence of a damping term, in which the conditions on the finite time extinction property fails (see Theorems~\ref{thm0s} and~\ref{thm0w} below). See Shimomura~\cite{MR2254620} for a related result with $m=1+\frac2N$.

We mention that it seems possible to apply the techniques of this paper to the consideration of some other complex-valued nonlinear equations such as the Gross-Pitaevskii equations, the Hartree-Fock equations, and the Ginzburg-Landau equations (see, e.g., B\'{e}gout and D\'{\i}az~\cite{MR3315701}, Antontsev, Dias and Figueira~\cite{MR3208711}, Okazawa and Yokota \cite{MR1886827} and its many references). 

Finally, we collect here some notations which will be used along with this paper. We let $\N_0=\N\cup\{0\}.$ Let $t\in\R.$ Then $t_+=\max\{t,0\}$ is the positive part of $t.$ We denote by $\ovl z$ the conjugate of the complex number $z,$ by $\Re(z)$ its real part and by $\Im(z)$ its imaginary part. For $1\le p\le\infty,$ $p^\prime$ is the conjugate of $p$ defined by $\frac{1}{p}+\frac{1}{p^\prime}=1.$ We write $\Gamma$ the boundary of a subset $\Omega\subset\R^N.$ Unless if specified, all functions are complex-valued $(H^1(\Omega)=H^1(\Omega;\C),$ etc). The notations $L^p(\Omega)$ $(p\in(0,\infty]),$ $W^{k,p}(\Omega),$ $W_0^{k,p}(\Omega),$ $H^k(\Omega),$ $H^k_0(\Omega)$ $(p\in[1,\infty],$ $k\in\N),$ $W^{-k,p^\p}(\Omega)$ and $H^{-k}(\Omega)$ $(p\in[1,\infty),$ $k\in\N)$ refer as the usual well known different Lebesgue, Sobolev and Hilbert spaces and their topological dual. By convention of notation, $W^{0,p}(\Omega)=W^{0,p}_0(\Omega)=L^p(\Omega).$ For a Banach space $X,$ we denote by $X^\star$ its topological dual and by $\langle\: . \; , \: . \:\rangle_{X^\star,X}\in\R$ the $X^\star-X$ duality product. In particular, for any $T\in L^{p^\prime}(\Omega)$ and $\vphi\in L^p(\Omega)$ with $1\le p<\infty,$ $\langle T,\vphi\rangle_{L^{p^\prime}(\Omega),L^p(\Omega)}=\Re\int_\Omega T(x)\ovl{\vphi(x)}\d x.$ The scalar product in $L^2(\Omega)$ between two functions $u,v$ is, $(u,v)_{L^2(\Omega)}=\Re\int_\Omega u(x)\ovl{v(x)}\d x.$ For a Banach space $X$ and $p\in[1,\infty],$ $u\in L^p_\loc\big([0,\infty);X\big)$ means that $u\in L^p_\loc\big((0,\infty);X\big)$ and for any $T>0,$ $u_{|(0,T)}\in L^p\big((0,T);X\big).$ In the same way, $u\in W^{1,p}_\loc\big([0,\infty);X\big)$ means that $u\in L^p_\loc\big([0,\infty);X\big),$ $u$ is absolutely continuous over $[0,\infty)$ (so it has a derivative $u^\p$ almost everywhere on $(0,\infty))$ and $u^\p\in L^p_\loc\big([0,\infty);X\big).$ For a real $x,$ $[x]$ denotes its integer part. As usual, we denote by $C$ auxiliary positive constants, and sometimes, for positive parameters $a_1,\ldots,a_n,$ write as $C(a_1,\ldots,a_n)$ to indicate that the constant $C$ depends only on $a_1,\ldots,a_n$ and that this dependence is continuous (we will use this convention for constants which are not denoted merely by ``$C$'').

\section{A semi-abstract result for finite time extinction}
\label{general}

We consider the following nonlinear Schrödinger equation.

\begin{empheq}[left=\empheqlbrace]{align}
	\label{nls}
	\vi\frac{\partial u}{\partial t}+\Delta u+a|u|^{-(1-m)}u=f(t,x),	&	\text{ in } (0,\infty)\times\Omega,			\\
	\label{nlsb}
	u(t)_{|\Gamma}=0,								&	\dfrac{}{} \text{ on } (0,\infty)\times\Gamma,	\\
	\label{u0}
	u(0)= u_0,										&	\text{ in } \Omega.
\end{empheq}
\medskip

\noindent
The next result proves the finite time extinction of solutions (in some cases even in the same time in which the source $f(t,x)$ vanishes) under suitable ``regularity'' conditions on the solution (this is the reason why we denote as ``semi-abstract'' such a framework). In the following sections we shall obtain sufficient conditions implying that such a framework holds.

\begin{thm}
\label{thmgenext}
Let $\Omega\subseteq\R^N$ be an open subset, $0<m\le1,$ $a\in\C,$ $f\in L^1_\loc\big([0,\infty);L^2(\Omega)\big)$ and $u_0\in L^2(\Omega).$ Assume that $u$ is any strong solution to~\eqref{nls}--\eqref{u0} $($see Definition~$\ref{defsol}$ below$)$ and that,
\begin{gather}
\label{thmgenext1}
u\in L^\infty\big((0,\infty);H^\ell_0(\Omega)\big),
\end{gather}
where $\ell=\left[\frac{N}2\right]+1$ $($or $H^\ell(\Omega)$ instead of $H^\ell_0(\Omega),$ if $\Omega$ is a half-space or if $\Omega$ has a bounded $C^{0,1}$-boundary$).$ Then the following conclusions hold.
\begin{enumerate}[$1)$]
\item
\label{pthmgenext1}
If there exists $T_0\ge0$ such that,
\begin{gather}
\label{thmgenext2}
\text{for almost every } t>T_0, \; f(t)=0,
\end{gather}
then there exists a finite time $T_\star\ge T_0$ such that,
\begin{gather}
\label{0}
\forall t\ge T_\star, \; \|u(t)\|_{L^2(\Omega)}=0.
\end{gather}
Furthermore,
\begin{gather}
\label{thmgenext2*}
T_\star\le\frac{2\,\ell\,C_\GN\,\|u\|_{L^\infty((0,\infty);H^\ell(\Omega))}^\frac{N(1-m)}{2\ell}}{\Im(a)(1-m)(2\ell-N)}\|u(T_0)\|_{L^2(\Omega)}^\frac{(1-m)(2\ell-N)}{2\ell}+T_0,
\end{gather}
where $C_\GN=C_\GN(N,m)$ is the constant in the inequality~\eqref{GN} below.
\item
\label{pthmgenext2}
There exist $\eps_\star=\eps_\star(\Im(a),N,m)$ satisfying the following property. Let $T_0>0$ and let $C_\GN$ be the constant in~\eqref{GN}. If,
\begin{gather}
\label{thmgenext3}
\|u\|_{L^\infty((0,\infty);H^\ell(\Omega))}^{1-m}\le\Im(a)\,C_\GN^{-1}\,\delta\,(1-\delta)\,T_0,
\end{gather}
and if for almost every $t>0,$
\begin{gather}
\label{thmgenext4}
\|f(t)\|_{L^2(\Omega)}^2\le\eps_\star\|u\|_{L^\infty((0,\infty);H^\ell(\Omega))}^{-\frac{2N}{2\ell-N}}\big(T_0-t\big)_+^\frac{2\delta-1}{1-\delta},
\end{gather}
where $\delta=\frac{(2\ell+N)+m(2\ell-N)}{4\ell}\in\left(\frac12,1\right),$ then~\eqref{0} holds true with $T_\star=T_0.$
\end{enumerate}
\end{thm}

\begin{rmk}
\label{rmkdelta}
Notice that $\delta\,(1-\delta)=\frac{(2\ell-N)(1-m)((2\ell+N)+m(2\ell-N))}{16\ell^2}$ and
$\frac{2\delta-1}{1-\delta}=2\frac{N(1-m)+2\ell m}{(2\ell-N)(1-m)}.$
\end{rmk}

\noindent
The following result collects several very useful \textit{a priori} estimates and some time differentiability conditions.

\begin{prop}
\label{propL2}
Let $\Omega\subseteq\R^N$ be an open subset, $0<m\le1,$ $a\in\C,$ $f\in L^1_\loc\big([0,\infty);L^2(\Omega)\big)$ and $u_0\in L^2(\Omega).$ Assume that $u$ is any weak solution to~\eqref{nls}--\eqref{u0} $($see Definition~$\ref{defsol}$ below$).$ Then we have the following results.
\begin{gather}
\label{Lm}
u\in L^{m+1}_\loc\big([0,\infty);L^{m+1}(\Omega)\big),	\\
\label{L2+}
\left\{
\begin{array}{rr}
\dfrac12\|u(t)\|_{L^2(\Omega)}^2+\Im(a)\dsp\vint_s^t\|u(\sigma)\|_{L^{m+1}(\Omega)}^{m+1}\d\sigma
												\ge\dfrac12\|u(s)\|_{L^2(\Omega)}^2		&			\\
+\,\Im\dsp\iint\limits_{s\,\Omega}^{\text{}\;\;t}f(\sigma,x)\,\ovl{u(\sigma,x)}\,\d x\,\d\sigma,	&	\text{if } \;	\Im(a)\le0,		\\
\dfrac12\|u(t)\|_{L^2(\Omega)}^2+\Im(a)\dsp\vint_s^t\|u(\sigma)\|_{L^{m+1}(\Omega)}^{m+1}\d\sigma
												\le\dfrac12\|u(s)\|_{L^2(\Omega)}^2		&			\\
+\,\Im\dsp\iint\limits_{s\,\Omega}^{\text{}\;\;t}f(\sigma,x)\,\ovl{u(\sigma,x)}\,\d x\,\d\sigma,	&	\text{if } \;	\Im(a)\ge0,
\end{array}
\right.
\end{gather}
for any $t\ge s\ge0.$ Finally, if $u$ satisfies one of the conditions below then the map $t\longmapsto\|u(t)\|_{L^2(\Omega)}^2$ belongs to $W^{1,1}_\loc\big([0,\infty);\R\big)$ and we have equality in~\eqref{L2+}.
\begin{enumerate}[$a)$]
\item
\label{ta}
$u$ is a strong solution $($see Definition~$\ref{defsol}$ below$),$
\item
\label{tb}
$|\Omega|<\infty,$
\item
\label{tc}
$m=1,$
\item
\label{td}
$\Im(a)=0.$
\end{enumerate}
\end{prop}

\begin{rmk}
\label{rmkthmgenext}
Here are some comments about Theorem~\ref{thmgenext}.
\begin{enumerate}[1)]
\item
\label{rmkthmgenext1}
Let $f$ satisfies~\eqref{thmgenext2} and let $u$ be a weak solution (see Definition~\ref{defsol} below). By~\eqref{L2+} we obtain that for any $t\ge T_0,$
\begin{gather*}
\begin{cases}
\|u(t)\|_{L^2(\Omega)}=\|u(T_0)\|_{L^2(\Omega)},	&	\text{if } \; \Im(a)=0,	\medskip \\
\|u(t)\|_{L^2(\Omega)}\ge\|u(T_0)\|_{L^2(\Omega)},	&	\text{if } \; \Im(a)<0.
\end{cases}
\end{gather*}
It follows that in those cases the finite time extinction is not reachable. If $m=1$ then we have, thanks to Proposition~\ref{propL2},
\begin{gather*}
\forall t\ge T_0, \; \|u(t)\|_{L^2(\Omega)}=\|u(T_0)\|_{L^2(\Omega)}e^{-\Im(a)(t-T_0)}.
\end{gather*}
And again, there is no finite time extinction.
\item
\label{rmkthmgenext2}
Let $u$ be a weak solution of \eqref{nls} (see Definition~\ref{defsol}). It is obvious from the equation and \ref{rmkthmgenext1}) of this remark that if $u$ vanishes at a finite time $T_\star\ge0$ then necessarily $f$ must satisfy \eqref{thmgenext2} (but not  necessarily the decay condition~\eqref{thmgenext4}) and that necessarily $\Im(a)>0$ and $m<1.$ If, in addition, $|\Omega|<\infty$ then we have,
\begin{gather}
\label{rmkthmgenext21}
T_\star\ge\frac{\|u(T_0\|_{L^2(\Omega)}^{1-m}}{(1-m)\Im(a)|\Omega|^\frac{1-m}2}+T_0.
\end{gather}
Indeed, it follows from \eqref{thmgenext2}, Proposition~\ref{propL2} and Hölder's inequality that for almost every $t>T_0,$
\begin{gather*}
\frac12\frac{\d}{\d t}\|u(t)\|_{L^2(\Omega)}^2=-\Im(a)\|u(t)\|_{L^{m+1}(\Omega)}^{m+1}\ge-\Im(a)|\Omega|^\frac{1-m}2\|u(t)\|_{L^2(\Omega)}^{m+1},
\end{gather*}
that is, $y^\p\ge-2\Im(a)|\Omega|^\frac{1-m}2y^\frac{m+1}2,$ where $y(\:.\:)=\|u(\:.\:)\|_{L^2(\Omega)}^2.$ After integration we get,
\begin{gather*}
y(t)^\frac{1-m}2\ge\left(y(T_0)^\frac{1-m}2-(1-m)\Im(a)|\Omega|^\frac{1-m}2(t-T_0)\right)_+,
\end{gather*}
for any $t\ge T_0,$ since $y\ge0.$ Hence the result.
\item
\label{rmkthmgenext3}
The proof of the finite time extinction of $u$ strongly relies on Gagliardo-Nirenberg's inequality (Lemma~\ref{lemGN} below), that is: for any $v\in H^\ell_0(\Omega)\cap L^{m+1}(\Omega)$ (or $H^\ell(\Omega)$ instead of $H^\ell_0(\Omega),$ if $\Omega$ is a half-space or if $\Omega$ has a bounded $C^{0,1}$-boundary),
\begin{gather}
\label{rmkthmgenext31}
\|v\|_{L^2(\Omega)}^\frac{(2\ell+N)+m(2\ell-N)}{2\ell}
\le C_\GN\|v\|_{L^{m+1}(\Omega)}^{m+1}\|v\|_{H^\ell(\Omega)}^\frac{N(1-m)}{2\ell},
\end{gather}
to get the ordinary differential inequality~\eqref{odi} below:
\begin{gather}
\label{rmkthmgenext32}
y^\p(t)+2\,\Im(a)\,C_\GN^{-1}\,\|u\|_{L^\infty((0,\infty);H^\ell(\Omega))}^{-\frac{N(1-m)}{2\ell}}\,y(t)^\delta\le0, \quad t>T_0,
\end{gather}
where $\delta=\frac{(2\ell+N)+m(2\ell-N)}{4\ell},$ $y=\|u(\:.\:)\|_{L^2(\Omega)}^2$ and $C_\GN=C_\GN(N,m,\ell).$ This holds thanks to the non increasing property~\eqref{L2+} of the mass (we recall that $\Im(a)>0$ is necessary to have finite time extinction, by \ref{rmkthmgenext1}) of this remark). But this method fails if $N\ge2\ell.$ Indeed, first of all, Gagliardo-Nirenberg's inequality imposes that $0\le m\le1.$ And as seen in \ref{rmkthmgenext1}) of this remark, finite time extinction is not reachable for $m=1.$ So, assume that $0\le m<1,$ \eqref{thmgenext2} is fulfilled and $u$ satisfies \eqref{thmgenext1}, where the integer $\ell$ has to be chosen later. Then for any $\ell\ge1,$ we may apply Lemma~\ref{lemGN} below, which is~\eqref{rmkthmgenext31} with $v=u(t),$ and we finally get \eqref{rmkthmgenext32}. But if $N$ is even and $\ell=\frac{N}2$ then $\delta=1$ and Lemma~\ref{lemodi} below yield,
\begin{gather}
\label{rmkthmgenext33}
\|u(t)\|_{L^2(\Omega)}\le\|u(T_0)\|_{L^2(\Omega)}e^{-\Im(a)\,C^{-1}\,(t-T_0)},
\end{gather}
for any $t\ge T_0,$ where $C=C(\|u\|_{L^\infty((0,\infty);H^\ell(\Omega))},N,m).$ In the same way, if $1\le\ell<\frac{N}2$ then $\delta>1$ and Lemma~\ref{lemodi} below yield,
\begin{gather}
\label{rmkthmgenext34}
\|u(t)\|_{L^2(\Omega)}\le\dfrac{\|u(T_0)\|_{L^2(\Omega)}}
{\left(1+\Im(a)\,C^{-1}(1-m)(N-2\ell)\|u(T_0)\|_{L^2(\Omega)}^\frac{(1-m)(N-2\ell)}{2\ell}(t-T_0)\right)^\frac{2\ell}{(1-m)(N-2\ell)}},
\end{gather}
for any $t\ge T_0,$ where $C=C(\|u\|_{L^\infty((0,\infty);H^\ell(\Omega))},N,m),$ and again this estimate does not give necessarily any finite time extinction result.
\end{enumerate}
\end{rmk}

\section{Finite time extinction and asymptotic behavior of solutions}
\label{finite}

Most of the results in this paper hold under the structural assumptions below.

\begin{assum}
\label{ass}
We assume that $\Omega\subseteq\R^N$ is a nonempty subset, $0<m\le1$ and $a\in\C$ with $\Im(a)>0.$ If $m<1$ then we assume further that,
\begin{gather}
\label{a}
2\sqrt m\,\Im(a)\ge(1-m)|\Re(a)|,	\\
\label{o}
|\Omega|<\infty.
\end{gather}
\end{assum}

\begin{thm}
\label{thmextH2}
Let Assumption~$\ref{ass}$ be fulfilled with $N\in\{1,2,3\}$ and $m<1.$ Let $f\in W^{1,1}_\loc\big([0,\infty);L^2(\Omega)\big),$ $u_0\in H^1_0(\Omega)$ and assume that one of the following hypotheses holds.
\begin{enumerate}[$1)$]
\item
\label{N1}
$N=1$ and $f\in W^{1,1}_\loc\big([0,\infty);H^1_0(\Omega)\big).$
\item
\label{N123}
$N\in\{1,2,3\},$ $\Omega$ is bounded with a $C^{1,1}$-boundary and $u_0\in H^2(\Omega)\cap H^1_0(\Omega).$ 
\end{enumerate}
Let $u$ be the unique strong solution of \eqref{nls}--\eqref{u0} $($see Definition~$\ref{defsol},$ Theorems~$\ref{thmstrongH1}$ and $\ref{thmstrongH2}$ and Remark~$\ref{rmkthmstrong}$ below$).$ Finally, assume that there exists $T_0\ge0$ such that,
\begin{gather*}
\text{for almost every } t>T_0, \; f(t)=0.
\end{gather*}
Then we have the following results.
\begin{enumerate}[$a)$]
\item
\label{thmextH21}
There exists a finite time $T_\star\ge T_0$ such that,
\begin{gather}
\label{thmext}
\forall t\ge T_\star, \; \|u(t)\|_{L^2(\Omega)}=0.
\end{gather}
Furthermore, $T_\star$ satisfies the estimates \eqref{thmgenext2*} and \eqref{rmkthmgenext21}.
\item
\label{thmextH22}
There exists $\eps_\star=\eps_\star(|a|,|\Omega|,N,m)$ satisfying the following property. Let $\delta$ be given in Property~$\ref{pthmgenext2})$ of Theorem~$\ref{thmgenext}.$ If $f\in W^{1,1}\big((0,\infty);H^1_0(\Omega)\big),$
\begin{gather*}
\begin{cases}
\left(\|u_0\|_{H^1_0(\Omega)}+\|f\|_{L^1((0,\infty);H^1_0(\Omega))}\right)^{1-m}\le\eps_\star\min\big\{1,T_0\big\},
&	\text{if } N=1,	\\
\left(\|u_0\|_{H^2(\Omega)}^m+\|f\|_{W^{1,1}((0,\infty);H^1_0(\Omega))}^m\right)^{1-m}\le\eps_\star\min\big\{1,T_0\big\},
&	\text{if } N\in\{2,3\},
\end{cases}
\end{gather*}
and if for almost every $t>0,$
\begin{gather*}
\|f(t)\|_{L^2(\Omega)}^2\le\eps_\star\big(T_0-t\big)_+^\frac{2\delta-1}{1-\delta},
\end{gather*}
then \eqref{thmext} holds with $T_\star=T_0.$
\end{enumerate}
\end{thm}

\begin{rmk}
\label{delta}
Notice that $\frac{2\delta-1}{1-\delta}=2\frac{1+m}{1-m},$ if $N\in\{1,2\}$ and $\frac{2\delta-1}{1-\delta}=2\frac{3+m}{1-m},$ if $N=3.$
\end{rmk}

\begin{rmk}
\label{rhmthmext}
Theorem~\ref{thmextH2} is an extension of the main result of Carles and Gallo~\cite{MR2765425} in the sense that they obtain the same conclusion as in $\ref{thmextH21})$ but under the additional conditions $\Re(a)=0,$ $f=0$ and without the lower bound for $T_\star.$ As far as we know, the result in $\ref{thmextH22})$ is new.
\end{rmk}

\noindent
The following result gives some asymptotic decay estimates, for large time, for the case of higher dimensions $N\ge4.$

\begin{thm}
\label{thm0s}
Let Assumption~$\ref{ass}$ be fulfilled with $N\ge4$ and $m<1.$ Let $f\in W^{1,1}_\loc\big([0,\infty);L^2(\Omega)\big)$ and let $u_0\in H^1_0(\Omega).$ Assume further that $f\in W^{1,1}_\loc\big([0,\infty);H^1_0(\Omega)\big)$ or $u_0\in H^2(\Omega)$ and that $\Omega$ is bounded with a $C^{1,1}$-boundary. Let $u$ be the unique strong solution of \eqref{nls}--\eqref{u0} $(see $ Definition~$\ref{defsol},$ Theorems~$\ref{thmstrongH1}$ and $\ref{thmstrongH2}$ and Remark~$\ref{rmkthmstrong}$ below$).$ Finally, assume that there exists $T_0\ge0$ such that
\begin{gather*}
\text{for almost every } t>T_0, \; f(t)=0.
\end{gather*}
Then we have for any $t\ge T_0,$
\begin{gather*}
\|u(t)\|_{L^2(\Omega)}\le\|u(T_0)\|_{L^2(\Omega)}e^{-\Im(a)\,C^{-1}\,(t-T_0)},
\end{gather*}
if $N=4$ and $u_0\in H^2(\Omega),$ and,
\begin{gather*}
\|u(t)\|_{L^2(\Omega)}\le\dfrac{\|u(T_0)\|_{L^2(\Omega)}}
{\left(1+\Im(a)\,C^{-1}(1-m)(N-2\ell)\|u(T_0)\|_{L^2(\Omega)}^\frac{(1-m)(N-2\ell)}{2\ell}(t-T_0)\right)^\frac{2\ell}{(1-m)(N-2\ell)}},
\end{gather*}
if $N\ge5$ or $u_0\in H^1_0(\Omega),$ where $C=C(\|u\|_{L^\infty((0,\infty);H^\ell(\Omega))},N,m).$
\end{thm}

\begin{thm}
\label{thm0w}
Let Assumption~$\ref{ass}$ be fulfilled, let $f\in L^1_\loc\big([0,\infty);L^2(\Omega)\big),$ let $u_0\in L^2(\Omega)$ and let $u$ be the unique weak solution of \eqref{nls}--\eqref{u0} $($see Definition~$\ref{defsol}$ and Theorem~$\ref{thmweak}$ below$).$ If
\begin{gather*}
f\in L^1\big((0,\infty);L^2(\Omega)\big),
\end{gather*}
then,
\begin{gather*}
\lim_{t\nearrow\infty}\|u(t)\|_{L^p(\Omega)}=0,
\end{gather*}
for any $p\in(0,2]$ $($with $p=2,$ if $m=1$ and $|\Omega|=\infty).$
\end{thm}

\begin{rmk}
\label{rmkthm0w}
Note that for $m=1$ in Theorem~\ref{thm0w}, if the stronger assumption~\eqref{thmgenext2} holds then we have,
\begin{gather*}
\forall t\ge T_0, \; \|u(t)\|_{L^2(\Omega)}=\|u(T_0)\|_{L^2(\Omega)}e^{-\Im(a)(t-T_0)}.
\end{gather*}
See~\ref{rmkthmgenext1}) of Remark~\ref{rmkthmgenext}.
\end{rmk}

\section{Existence and uniqueness of solutions}
\label{exiuni}

Here and after, we shall always identify $L^2(\Omega)$ with its topological dual. Let $\Omega\subseteq\R^N$ be an open subset, let $0<m\le1$ and let $X=H\cap L^{m+1}(\Omega),$ where $H=L^2(\Omega)$ or $H=H^1_0(\Omega).$ It follows from Lemma~\ref{lemden0} and \ref{embL21}) of Lemma~\ref{lemC} below that,
\begin{gather*}
X^\star=H^\star+L^\frac{m+1}m(\Omega),		\\
L^{m+1}_\loc\big([0,\infty);X\big)\cap W^{1,\frac{m+1}{m}}_\loc\big([0,\infty);X^\star\big)\inj C\big([0,\infty);L^2(\Omega)\big).
\end{gather*}
This justifies the notion of solution below (and it explains the sense in which the initial condition is satisfied).

\begin{defi}
\label{defsol}
Let $\Omega\subseteq\R^N$ be an open subset, $0<m\le1,$ $a\in\C,$ $f\in L^1_\loc\big([0,\infty);L^2(\Omega)\big)$ and $u_0\in L^2(\Omega).$ Let us consider the following assertions.
\begin{enumerate}[1)]
\item
\label{defsol1}
$u\in L^{m+1}_\loc\big([0,\infty);H^1_0(\Omega)\cap L^{m+1}(\Omega)\big)\cap W^{1,\frac{m+1}{m}}_\loc\big([0,\infty);H^\star+L^\frac{m+1}m(\Omega)\big).$
\item
\label{defsol2}
For almost every $t>0,$ $\Delta u(t)\in H^\star.$
\item
\label{defsol3}
$u$ satisfies~\eqref{nls} in $\Dr^\p\big((0,\infty)\times\Omega\big).$
\item
\label{defsol4}
$u(0)=u_0.$
\end{enumerate}
We shall say that $u$ is a \textit{strong solution} if $u$ is a $H^2$-solution or a $H^1_0$-solution. We shall say that $u$ is a $H^2$-\textit{solution of} \eqref{nls}--\eqref{u0} \big(respectively, a $H^1_0$-\textit{solution of} \eqref{nls}--\eqref{u0}\big), if $u$ satisfies the Assertions~\ref{defsol1})--\ref{defsol4}) with $H=L^2(\Omega)$ \big(respectively, with $H=H^1_0(\Omega)\big).$
\\
We shall say that $u$ is a $L^2$-\textit{solution} or simply a \textit{weak solution} of \eqref{nls}--\eqref{u0} is there exists a pair,
\begin{gather}
\label{fn}
(f_n,u_n)_{n\in\N}\subset L^1_\loc\big([0,\infty);L^2(\Omega)\big)\times C\big([0,\infty);L^2(\Omega)\big),
\end{gather}
such that for any $n\in\N,$ $u_n$ is a $H^2$-solution of \eqref{nls}--\eqref{nlsb} where the right-hand side member of \eqref{nls} is $f_n,$ and if
\begin{gather}
\label{cv}
f_n\xrightarrow[n\to\infty]{L^1((0,T);L^2(\Omega))}f \; \text{ and } \; u_n\xrightarrow[n\to\infty]{C([0,T];L^2(\Omega))}u,
\end{gather}
for any $T>0.$
\end{defi}

\begin{rmk}
\label{rmkdefsol}
Before making some comments on the above definition, it is useful to analyze some peculiar properties which arise when $\Omega$ is unbounded. Let $0<m\le1.$ Set for any $z\in\C,$ $g(z)=|z|^{-(1-m)}z$ $(g(0)=0)$ and let us define the mapping for any measurable function $u:\Omega\tends\C,$ which we still denote by $g,$ by $g(u)(x)=g(u(x)).$ Let $H=L^2(\Omega)$ or $H=H^1_0(\Omega).$ It follows from~\eqref{lemmon2} below that,
\begin{gather}
\label{gm}
g\in C\big(L^{m+1}(\Omega);L^\frac{m+1}m(\Omega)\big) \text{ and } g \text{ is bounded on bounded sets.}
\end{gather}
In particular, if $|\Omega|<\infty$ or if $m=1$ then $H^1_0(\Omega)\inj L^2(\Omega)\inj L^{m+1}(\Omega)$ with dense embedding and thus, $L^\frac{m+1}m(\Omega)\inj L^2(\Omega)\inj H^{-1}(\Omega).$ We then obtain,
\begin{gather}
\label{g2}
g\in C\big(L^2(\Omega);L^2(\Omega)\big)\cap C\big(H^1_0(\Omega);H^{-1}(\Omega)\big)
 \text{ and } g \text{ is bounded on bounded sets,}
\end{gather}
and Assertion~\ref{defsol1}) becomes,
\begin{gather}
\label{defsol1bis}
u\in L^{m+1}_\loc\big([0,\infty);H^1_0(\Omega)\big)\cap W^{1,\frac{m+1}{m}}_\loc\big([0,\infty);H^\star\big).
\end{gather}
But if $|\Omega|=\infty$ and $m<1$ then the regularity~\eqref{g2} is not anymore valid. By Lemma~\ref{lemden0} below, we have,
\begin{gather}
\label{den1}
\Dr(\Omega)\inj X\inj L^{m+1}(\Omega) \text{ with both dense embeddings},
\end{gather}
where $X=H\cap L^{m+1}(\Omega).$ It follows that,
\begin{gather}
\label{den2}
L^\frac{m+1}m(\Omega)\inj X^\star\inj\Dr^\p(\Omega).
\end{gather}
This gives with~\eqref{gm},
\begin{gather}
\label{gmm}
g\in C(X,X^\star) \text{ and } g \text{ is bounded on bounded sets.}
\end{gather}
It follows from \eqref{gm} and \eqref{den1}--\eqref{gmm} that,
\begin{gather}
\label{dualg}
\langle g(u),v\rangle_{X^\star,X}=\langle g(u),v\rangle_{L^\frac{m+1}m(\Omega),L^{m+1}(\Omega)}
=\Re\vint_\Omega g(u)\ovl v\d x,
\end{gather}
for any $u,v\in X.$ Now, let us make some comments about Definition~\ref{defsol}. 
\begin{enumerate}[1)]
\item
\label{rmkdefsol1}
As seen at the beginning of this section, any strong or weak solution belongs to $C\big([0,\infty);L^2(\Omega)\big)$ and Assertion~\ref{defsol4}) makes sense in $L^2(\Omega).$
\item
\label{rmkdefsol2}
It is obvious that a $H^2$-solution is also a $H^1_0$-solution and a weak solution. But it is not clear that a $H^1_0$-solution is a weak solution, without assuming a continuous dependence of the solution with respect to the initial data. Such a result will be established with the additional assumption~\eqref{a} on $a$ (see Lemma~\ref{lemdep} below).
\item
\label{rmkdefsol3}
If $|\Omega|<\infty$ or if $m=1$ then it follows from~\eqref{g2}, \eqref{defsol1bis} and Assertion~\ref{defsol2}) that any $H^2$-solution (respectively, any $H^1_0$-solution) satisfies~\eqref{nls} in $L^2(\Omega)$ \big(respectively, in $H^{-1}(\Omega)\big),$ for almost every $t>0.$ Note also that Assertion~\ref{defsol2}) of Definition~\ref{defsol} is not an additional assumption for the $H^1_0$-solutions.
\item
\label{rmkdefsol4}
If $|\Omega|=\infty$ and if $m<1$ then it follows from~\eqref{gmm} and Assertions~\ref{defsol1}) and \ref{defsol2}) that any $H^2$-solution (respectively, any $H^1_0$-solution) satisfies~\eqref{nls} in $L^2(\Omega)+L^\frac{m+1}m(\Omega)$ \big(respectively, in $H^{-1}(\Omega)+L^\frac{m+1}m(\Omega)\big),$ for almost every $t>0.$
\item
\label{rmkdefsol5}
Assume that $u$ is a weak solution. By Definition~\ref{defsol}, there exists $(f_n,u_n)_{n\in\N}$ satisfying~\eqref{fn}--\eqref{cv} such that for any $n\in\N,$ $u_n$ is a $H^2$-solution of \eqref{nls}--\eqref{nlsb} where the right-hand side of \eqref{nls} is $f_n.$ Applying \eqref{lemmon2}--\eqref{lemmon3} below, we deduce that for any $T>0,$
\begin{align*}
&	\Delta u_n\xrightarrow[n\to\infty]{C([0,T];H^{-2}(\Omega))}\Delta u,					\\
&	g(u_n)\xrightarrow[n\to\infty]{C([0,T];L^2(\Omega))}g(u), \; \text{ if } \; |\Omega|<\infty,	\\
&	g(u_n)\xrightarrow[n\to\infty]{C([0,T];L^\frac2m(\Omega))}g(u).
\end{align*}
Now, we set: $Y=H^2_0(\Omega)\cap L^\frac2{2-m}(\Omega).$ By Lemma~\ref{lemden0} below, we have,
\begin{gather*}
Y^\star=H^{-2}(\Omega)+L^\frac2m(\Omega),													\\
\Dr(\Omega)\inj Y\inj H^2_0(\Omega),L^2(\Omega),L^\frac2{2-m}(\Omega) \; \text{ with dense embedding},	\\
H^{-2}(\Omega),L^2(\Omega),L^\frac2m(\Omega)\inj Y^\star\inj\Dr^\p(\Omega).
\end{gather*}
Using the above uniform convergences and~\eqref{cv}, we deduce that,
\begin{gather*}
\vint_0^\infty\left\langle\vi\frac{\partial u}{\partial t}+\Delta u
+ag(u),\vphi\right\rangle_{Y^\star,Y}\psi(t)\,\d t
=\vint_0^\infty\big\langle f(t),\vphi\big\rangle_{Y^\star,Y}\psi(t)\d t.
\end{gather*}
for any $\vphi\in Y$ and $\psi\in\Dr\big((0,\infty);\R\big).$ \\
As a conclusion, if $u$ is a weak solution then $u\in W^{1,1}_\loc\big([0,\infty);Y^\star\big)$ and it solves~\eqref{nls} in $Y^\star,$ for almost every $t>0.$ In particular, $u$ satisfies~\eqref{nls} in $\Dr^\p\big((0,\infty)\times\Omega\big).$ If, in addition, $|\Omega|<\infty$ or if $m=1$ then we deduce from the above that $u\in W^{1,1}_\loc\big([0,\infty);H^{-2}(\Omega)\big)$ and $u$ solves~\eqref{nls} in $H^{-2}(\Omega),$ for almost every $t>0.$
\item
\label{rmkdefsol6}
When $m<1$ then except for Theorem~\ref{thmgenext} and Proposition~\ref{propL2}, all the results of the following Sections~\ref{general}--\ref{exiuni} will be stated with $|\Omega|<\infty.$
\item
\label{rmkdefsol7}
Notice that the boundary condition $u(t)_{|\Gamma}=0$ is included in the assumption $u(t)\in H^1_0(\Omega).$ 
\end{enumerate}
\medskip
\end{rmk}

\begin{thm}[\textbf{Existence and uniqueness of $\bs{L^2}$-solutions}]
\label{thmweak}
Let Assumption~$\ref{ass}$ be fulfilled and let $f\in L^1_\loc\big([0,\infty);L^2(\Omega)\big).$ Then for any $u_0\in L^2(\Omega),$ there exists a unique weak solution $u$ to \eqref{nls}--\eqref{u0}. In addition, we have the following properties.
\begin{enumerate}[$1)$]
\item
\label{thmweak1}
The map $t\longmapsto\|u(t)\|_{L^2(\Omega)}^2$ belongs to $W^{1,1}_\loc\big([0,\infty);\R\big)$ and we have,
\begin{gather}
\label{L2}
\frac12\frac{\d}{\d t}\|u(t)\|_{L^2(\Omega)}^2+\Im(a)\|u(t)\|_{L^{m+1}(\Omega)}^{m+1}=\Im\vint_\Omega f(t,x)\,\ovl{u(t,x)}\,\d x,
\end{gather}
for almost every $t>0.$
\item
\label{thmweak2}
If $v$ is another weak solution of \eqref{nls}--\eqref{nlsb} with $v(0)=v_0\in L^2(\Omega)$ and $h\in L^1_\loc([0,\infty);L^2(\Omega)),$ instead of $f$ in \eqref{nls} then,
\begin{gather}
\label{estthmweak}
\|u(t)-v(t)\|_{L^2(\Omega)}\le\|u(s)-v(s)\|_{L^2(\Omega)}+\vint_s^t\|f(\sigma)-h(\sigma)\|_{L^2(\Omega)}\d\sigma,
\end{gather}
for any $t\ge s\ge0.$
\end{enumerate}
\end{thm}

\begin{thm}[\textbf{Existence and uniqueness of $\bs{H^1_0}$-solutions}]
\label{thmstrongH1}
Let Assumption~$\ref{ass}$ be fulfilled and let $f\in W^{1,1}_\loc\big([0,\infty);H^1_0(\Omega)\big).$ Then for any $u_0\in H^1_0(\Omega),$ there exists a unique $H^1_0$-solution $u$ to \eqref{nls}--\eqref{u0}. Furthermore, $u$ is also a weak solution and satisfies the following properties.
\begin{enumerate}[$1)$]
\item
\label{thmstrongH11}
$u\in C\big([0,\infty);L^2(\Omega)\big)\cap C^1\big([0,\infty);H^{-2}(\Omega)\big)$ and $u$ satisfies \eqref{nls} in $H^{-2}(\Omega),$ for any $t\ge0.$
\item
\label{thmstrongH12}
$u\in C_\w\big([0,\infty);H^1_0(\Omega)\big)\cap W^{1,\infty}_\loc\big([0,\infty);H^{-1}(\Omega)\big)$ and,
\begin{empheq}[left=\empheqlbrace]{align}
\label{strongH11}
&	\|u(t)-u(s)\|_{L^2(\Omega)}\le M|t-s|^\frac12,	\\
\label{strongH12}
&	\|\nabla u(t)\|_{L^2(\Omega)}\le\|\nabla u_0\|_{L^2(\Omega)}+\int_0^t\|\nabla f(s)\|_{L^2(\Omega)}\d s,
\end{empheq}
for any $t\ge s\ge0,$ where $M^2=2\|u\|_{L^\infty((s,t);H^1_0(\Omega))}\|u_t\|_{L^\infty((s,t);H^{-1}(\Omega))}.$
\item
\label{thmstrongH13}
The map $t\longmapsto\|u(t)\|_{L^2(\Omega)}^2$ belongs to $C^1\big([0,\infty);\R\big)$ and~\eqref{L2} holds for any $t\ge0.$
\item
\label{thmstrongH14}
If $f\in W^{1,1}\big((0,\infty);H^1_0(\Omega)\big)$ then we have,
\begin{gather*}
u\in L^\infty\big((0,\infty);H^1_0(\Omega)\big)\cap W^{1,\infty}\big((0,\infty);H^{-1}(\Omega)\big)
\cap C^1_\b\big([0,\infty);H^{-2}(\Omega)\big).
\end{gather*}
\end{enumerate}
\end{thm}

\begin{thm}[\textbf{Existence and uniqueness of $\bs{H^2}$-solutions}]
\label{thmstrongH2}
Let Assumption~$\ref{ass}$ be fulfilled and let $f\in W^{1,1}_\loc\big([0,\infty);L^2(\Omega)\big).$ Then for any $u_0\in H^1_0(\Omega)$ with $\Delta u_0\in L^2(\Omega),$ there exists a unique $H^2$-solution $u$ to \eqref{nls}--\eqref{u0}. Furthermore, $u$ satisfies the following properties.
\begin{enumerate}[$1)$]
\item
\label{thmstrongH21}
$u\in C\big([0,\infty);H^1_0(\Omega)\big)\cap C^1\big([0,\infty);H^{-1}(\Omega)\big),$ $u$ satisfies \eqref{nls} in $H^{-1}(\Omega),$ for any $t\ge0.$
\item
\label{thmstrongH22}
$u\in W^{1,\infty}_\loc\big([0,\infty);L^2(\Omega)\big),$ $\Delta u\in L^\infty_\loc\big([0,\infty);L^2(\Omega)\big)$ and,
\begin{empheq}[left=\empheqlbrace]{align}
\label{strongH21}
	&	\|u(t)-u(s)\|_{L^2(\Omega)}\le\|u_t\|_{L^\infty((s,t);L^2(\Omega))}|t-s|,			\dfrac{}{}				\\
\label{strongH22}
	&	\|\nabla u(t)-\nabla u(s)\|_{L^2(\Omega)}\le M|t-s|^\frac12,										\\
\label{strongH23}
	&	\left\|u_t\right\|_{L^\infty((0,t);L^2(\Omega))}
		\le\|\Delta u_0+a|u_0|^{m-1}u_0-f(0)\|_{L^2(\Omega)}+\int_0^t\|f^\p(\sigma)\|_{L^2(\Omega)}\d\sigma,
\end{empheq}
for any $t\ge s\ge0,$ where $M^2=2\|u_t\|_{L^\infty((s,t);L^2(\Omega))}\|\Delta u\|_{L^\infty((s,t);L^2(\Omega))}.$
\item
\label{thmstrongH23}
The map $t\longmapsto\|u(t)\|_{L^2(\Omega)}^2$ belongs to $C^1\big([0,\infty);\R\big)$ and~\eqref{L2} holds for any $t\ge0.$
\item
\label{thmstrongH24}
If $f\in W^{1,1}\big((0,\infty);L^2(\Omega)\big)$ then we have,
\begin{align*}
	&	u\in C_\b\big([0,\infty);H^1_0(\Omega)\big)\cap C_\b^1\big([0,\infty);H^{-1}(\Omega)\big)\cap
											W^{1,\infty}\big((0,\infty);L^2(\Omega)\big),		\\
	&	\Delta u\in L^\infty\big((0,\infty);L^2(\Omega)\big).
\end{align*}
\end{enumerate}
\end{thm}

\begin{rmk}
\label{rmkthmstrong}
Let $E=\left\{u\in H^1_0(\Omega); \Delta u\in L^2(\Omega)\right\}$ with $\|u\|_E^2=\|u\|_{L^2(\Omega)}^2+\|\Delta u\|_{L^2(\Omega)}^2.$ We recall that $E\subset H_\loc^2(\Omega)$ (Theorem~8.8, p.183-184, in Gilbarg and Trudinger~\cite{MR1814364}). If $\Omega=\R^N$ then $E=H^2(\R^N)$ with equivalent norms (by the Fourier transform and Plancherel's formula), while if $\Omega$ is bounded and $\Gamma$ is of class $C^{1,1}$ then $E=H^2(\Omega)\cap H^1_0(\Omega)$ with equivalent norms (Theorem~8.12, p.186, in Gilbarg and Trudinger~\cite{MR1814364} and Corollary 2.5.2.2, p.131, in Grisvard~\cite{MR3396210}). Note that for the equivalence of the norms, we may use the inequalities,
\begin{gather}
\label{E}
\|\nabla u\|_{L^2(\Omega)}^2\le\|u\|_{L^2(\Omega)}\|\Delta u\|_{L^2(\Omega)}
\le\|u\|_{L^2(\Omega)}^2+\|\Delta u\|_{L^2(\Omega)}^2,
\end{gather}
which hold for any subset $\Omega\subseteq\R^N$ and any $u\in H^2(\Omega)\cap H^1_0(\Omega).$
\end{rmk}

\begin{rmk}
\label{rmkf0}
Since $f\in C\big([0,\infty);L^2(\Omega)\big)$ (by \ref{W11}) of Lemma~\ref{lemC}), estimate~\eqref{strongH23} with $f(0)$ makes sense.
\end{rmk}

\begin{rmk}
\label{rmkthmstrongH10}
It follows from~\eqref{estthmweak} and \eqref{strongH12} that if $N=1$ then the decay assumptions~\eqref{thmgenext3} and \eqref{thmgenext4} may be replaced with,
\begin{gather}
\nonumber
\left(\|u_0\|_{H^1_0(\Omega)}+\|f\|_{L^1((0,\infty);H^1_0(\Omega))}\right)^{1-m}\le\eps_\star\min\big\{1,T_\star\big\},	\\
\label{rmkthmstrongH101}
\|f(t)\|_{L^2(\Omega)}^2\le\eps_\star\big(T_\star-t\big)_+^\frac{2\delta-1}{1-\delta},
\end{gather}
for almost every $t>0,$ where $\eps_\star=\eps_\star(\Im(a),N,m).$ In the same way, it follows from~\eqref{estthmweak}, \eqref{strongH12}, \eqref{strongH23}, Remark~\ref{rmkthmstrong} and~\eqref{nls} that if $N\le3$ and $\Omega$ is bounded with a $C^{1,1}$-boundary then~\eqref{thmgenext3} may be replaced with,
\begin{gather*}
\left(\|u_0\|_{H^2(\Omega)}^m+\|f\|_{W^{1,1}((0,\infty);H^1_0(\Omega))}^m\right)^{1-m}\le\eps_\star\min\big\{1,T_\star\big\},
\end{gather*}
and~\eqref{thmgenext4} with \eqref{rmkthmstrongH101}, where $\eps_\star=\eps_\star(|a|,|\Omega|,N,m).$
\end{rmk}

\section{Proof of the semi-abstract result on the finite time extinction}
\label{prooffinite}

\noindent
The proof of Theorem~\ref{thmgenext} relies on the three following lemmas.

\begin{lem}
\label{lemodi}
Let $y\in W^{1,1}_\loc\big([0,\infty);\R\big)$ with $y\ge0$ over $(0,\infty),$ $\delta\in\R,$ $\alpha>0$ and $T_0\ge0.$ If
\begin{gather*}
y^\p+2\alpha y^\delta\le0,
\end{gather*}
almost everywhere on $(T_0,\infty),$ then we have,
\begin{gather*}
y(t)\le
\begin{cases}
\Big(y(T_0)^{1-\delta}+2\alpha(1-\delta)(T_0-t)\Big)^\frac1{1-\delta}_+,			&	\text{if } \delta<1,	\medskip \\
y(T_0)e^{-2\alpha(t-T_0)},												&	\text{if } \delta=1,	\medskip \\
\dfrac{y(T_0)}{\big(1+2\alpha(\delta-1)y(T_0)^{\delta-1}(t-T_0)\big)^\frac1{\delta-1}}	&	\text{if } \delta>1,
\end{cases}
\end{gather*}
for any $t\ge T_0.$ In particular, if $\delta<1$ then for any $t\ge T_\star,$ $y(t)=0$ where,
\begin{gather*}
T_\star\le\frac1{2\alpha(1-\delta)}y(T_0)^{1-\delta}+T_0.
\end{gather*}
\end{lem}

\begin{proof*}
The result follows by integration of the ordinary differential inequality over $(T_0,t).$
\medskip
\end{proof*}

\noindent
The following lemma improves a similar result contained in Antontsev, D{\'{\i}}az and Shmarev~\cite{MR2002i:35001} (Proposition~1.1, p.77, and its proof, p.75--77).

\begin{lem}
\label{lemodif}
Let $y\in W^{1,1}_\loc\big([0,\infty);\R\big)$ with $y\ge0$ over $[0,\infty),$ $\delta\in(0,1),$ $\alpha,T_0>0$ and,
\begin{gather}
\label{lemodif1}
y_\star=\left(\alpha\,\delta^\delta(1-\delta)\right)^\frac1{1-\delta},		\\
\label{lemodif2}
x_\star=\left(\alpha\,\delta\,(1-\delta)\,T_0\right)^\frac1{1-\delta}.
\end{gather}
If,
\begin{gather}
\label{lemodif3}
y(0)\le x_\star,
\end{gather}
and if for almost every $t>0,$
\begin{gather}
\label{lemodif4}
y^\p(t)+\alpha y(t)^\delta\le y_\star\left(T_0-t\right)_+^\frac\delta{1-\delta},
\end{gather}
then for any $t\ge T_0,$ $y(t)=0.$
\end{lem}

\begin{proof*}
Set for any $t\in[0,T_0],$ $z(t)=x_\star T_0^{-\frac1{1-\delta}}\left(T_0-t\right)^\frac1{1-\delta}.$ We have for almost every $t\in(0,T_0),$
\begin{gather}
\label{demlemodif}
z^\p(t)+\alpha z(t)^\delta=y_\star\left(T_0-t\right)^\frac\delta{1-\delta}\ge y^\p(t)+\alpha y(t)^\delta.
\end{gather}
We claim that for any $t\in[0,T_0],$ $y(t)\le z(t).$ If not, since by \eqref{lemodif3} $z(0)\ge y(0)$ and $y$ and $z$ are continuous over $[0,T_0]$ (by \ref{W11}) of Lemma~\ref{lemC}), there exist $t_\star\in[0,T_0)$ and $\eps\in(0,T_0-t_\star)$ such that $y(t_\star)=z(t_\star)$ and $y(t)>z(t),$ for any $t\in(t_\star,t_\star+\eps).$ This leads with \eqref{demlemodif} to, $y^\p\le z^\p,$ almost everywhere on $(t_\star,t_\star+\eps).$ Integrating over $(t_\star,t)$ for $t\in(t_\star,t_\star+\eps),$ we obtain that $y(t)\le z(t),$ for any $t\in[t_\star,t_\star+\eps].$ A contradiction. Hence the claim. In particular, $y(T_0)\le z(T_0)=0.$ But from \eqref{lemodif4}, $y$ is non increasing over $(T_0,\infty).$ Hence the result, since $y\ge0$ everywhere.
\medskip
\end{proof*}

\begin{rmk}
\label{rmklemodif}
Let us explain how we found $y_\star$ and $x_\star$ in Lemma~\ref{lemodif}. We look for a solution of the ordinary differential inequality \eqref{lemodif4}. Set for any $x\ge0,$
\begin{gather*}
\forall x\ge0, \;
f(x)=(1-\delta)^{-1}T_0^{-\frac1{1-\delta}}x^\delta\left(\alpha(1-\delta)T_0-x^{1-\delta}\right),	\\
\forall t\in[0,T_0], \; z(t)=xT_0^{-\frac1{1-\delta}}\left(T_0-t\right)_+^\frac1{1-\delta}.
\end{gather*}
We want $z(0)=x\ge y(0)$ to apply our proof. A straightforward calculation yields,
\begin{gather*}
z^\p(t)+\alpha z(t)^\delta=f(x)\left(T_0-t\right)^\frac\delta{1-\delta}.
\end{gather*}
We compute, $\underset{x\ge0}{\argmax}\,f(x)=x_\star,$ where $x_\star$ is given by \eqref{lemodif2}, and $f(x_\star)=y_\star,$ where $y_\star$ is given by \eqref{lemodif1}. We then choose $x=x_\star$ in the definition of $z$ and we obtain the condition \eqref{lemodif3}.
\end{rmk}

\begin{lem}[\textbf{Gagliardo-Nirenberg's inequality}]
\label{lemGN}
Let $N\in\N,$ let $\Omega\subseteq\R^N$ be an open subset, let $0\le m\le1$ and let $\ell\in\N.$ Then for any $v\in H^\ell_0(\Omega)\cap L^{m+1}(\Omega),$
\begin{gather}
\label{GN}
\|v\|_{L^2(\Omega)}^\frac{(2\ell+N)+m(2\ell-N)}{2\ell}
\le C\|v\|_{L^{m+1}(\Omega)}^{m+1}\|v\|_{H^\ell(\Omega)}^\frac{N(1-m)}{2\ell},
\end{gather}
where $C=C(m,\ell,N).$ If $\Omega$ is a half-space or if $\Omega$ has a bounded $C^{0,1}$-boundary then \eqref{GN} holds for any $v\in H^\ell(\Omega).$
\end{lem}

\begin{proof*}
See, for instance, Friedman~\cite{MR0445088}, Theorem~9.3, p.24, for $v\in\Dr(\R^N)$ and so, by extension and density, for $v\in H^\ell_0(\Omega)\cap L^{m+1}(\Omega).$ If $\Omega$ is a half-space or if $\Omega$ has a bounded $C^{0,1}$-boundary then there exists a linear extension operator $E$ such that for any $k\in\N_0$ and $p\in[1,\infty],$
\begin{gather*}
E\in\Lr\big(W^{k,p}(\Omega);W^{k,p}(\R^N)\big),
\end{gather*}
and $Eu=u,$ almost everywhere in $\Omega$ (Stein~\cite{MR0290095}, Theorem~5 and \S3.2, p.181 and \S 3.3, p.189; Adams~\cite{MR0450957}, Theorem~4.26, p.84; see also Grisvard~\cite{MR3396210}, Theorem~1.4.3.1, p.25).
\medskip
\end{proof*}

\begin{vproof}{of Proposition~\ref{propL2}.}
Let the assumptions of the theorem be fulfilled. We first assume that $u$ is a strong solution. Let $H$ be as in Definition~\ref{defsol} and let $X=H\cap L^{m+1}(\Omega).$ By Definition~\ref{defsol}, we have \eqref{Lm} and by \ref{rmkdefsol3}) and \ref{rmkdefsol4}) of Remark~\ref{rmkdefsol}, we can take the $X^\star-X$ duality product with $\vi u.$ Estimate~\eqref{L2+} with equality then follows from \eqref{dualg} and \ref{lemmass}) of Lemma~\ref{lemmassene}. Now, assume that $u$ is a weak solution. Let $(f_n)_{n\in\N}$ and $(u_n)_{n\in\N}$ be as in Definition~\ref{defsol}. According to the above, it follows from Hölder's inequality that $f\ovl u\in L^1_\loc\big([0,\infty);L^1(\Omega)\big)$ and,
\begin{gather}
\label{demthmgenext1}
f_n\ovl{u_n}\xrightarrow[n\to\infty]{L^1_\loc([0,\infty);L^1(\Omega))}f\ovl u,\\
\begin{split}
\label{demthmgenext2}
	&	\; \frac12\|u_n(t)\|_{L^2(\Omega)}^2+\Im(a)\vint_s^t\|u_n(\sigma)\|_{L^{m+1}(\Omega)}^{m+1}\d\sigma		\\
   =	&	\; \frac12\|u_n(s)\|_{L^2(\Omega)}^2+\Im\iint\limits_{s\,\Omega}^{\text{}\;\;t}f_n(\sigma,x)\,\ovl{u_n(\sigma,x)}\,\d x\,\d\sigma,
\end{split}
\end{gather}
for any $n\in\N$ and $t\ge s\ge0.$ If $|\Omega|<\infty$ or if $m=1$ then for any $T>0,$ $C([0,T];L^2(\Omega))\inj C([0,T];L^{m+1}(\Omega))$ and then we are allowed to pass to the limit in \eqref{demthmgenext2} under the integral symbol. We then get with \eqref{demthmgenext1} the desired result under the hypotheses \ref{tb}), \ref{tc}) or \ref{td}). If $|\Omega|=\infty,$ $m<1$ and $\Im(a)\ge0$ then for any $T>0,$ $C([0,T];L^2(\Omega))\inj C([0,T];L^{m+1}_\loc(\Omega)).$ By \eqref{demthmgenext2},
\begin{multline*}
\frac12\|u_n(t)\|_{L^2(\Omega)}^2+\Im(a)\vint_s^t\|u_n(\sigma)\|_{L^{m+1}(\Omega\cap B(0,R))}^{m+1}\d\sigma			\\
\le\frac12\|u_n(s)\|_{L^2(\Omega)}^2+\Im\iint\limits_{s\,\Omega}^{\text{}\;\;t}f_n(\sigma,x)\,\ovl{u_n(\sigma,x)}\,\d x\,\d\sigma,
\end{multline*}
for any $t>s>0,$ $R>0$ and $n\in\N.$ Passing to the limit in $n$ first and then in $R$ then, we obtain \eqref{Lm} and \eqref{L2+} with the help of the monotone convergence Theorem and \eqref{demthmgenext1}. We proceed in the same way if $|\Omega|=\infty,$ $m<1$ and $\Im(a)\le0.$
\medskip
\end{vproof}

\begin{vproof}{of Theorem~\ref{thmgenext}.}
By \eqref{GN} and Proposition~\ref{propL2}, we have for almost every $t>0,$
\begin{gather*}
\|u(t)\|_{L^2(\Omega)}^\frac{(2\ell+N)+m(2\ell-N)}{2\ell}
\le C_\GN\|u\|_{L^\infty((0,\infty);H^\ell(\Omega))}^\frac{N(1-m)}{2\ell}\|u(t)\|_{L^{m+1}(\Omega)}^{m+1},	\\
\dfrac{\d}{\d t}\|u(t)\|_{L^2(\Omega)}^2+2\Im(a)\|u(t)\|_{L^{m+1}(\Omega)}^{m+1}=2\Im\vint_\Omega f(t,x)\ovl{u(t,x)}\d x.
\end{gather*}
It follows that,
\begin{gather}
\label{demthmgenext3}
\dfrac{\d}{\d t}\|u(t)\|_{L^2(\Omega)}^2+2\alpha\|u(t)\|_{L^2(\Omega)}^{2\delta}
\le2\vint_\Omega|f(t,x)||u(t,x)|\d x,
\end{gather}
for almost every $t>0,$ where $\alpha=\Im(a)C_\GN^{-1}\|u\|_{L^\infty((0,\infty);H^\ell(\Omega))}^{-\frac{N(1-m)}{2\ell}}$ and $\delta=\frac{(2\ell+N)+m(2\ell-N)}{4\ell}.$ Since $0<m<1$ and $\ell=\left[\frac{N}2\right]+1,$ we have $\frac12<\delta<1.$ Using the Young inequality,
\begin{gather*}
xy\le\frac{\eps^{-p^\p}}{p^\p}x^{p^\p}+\frac{\eps^p}py^p,
\end{gather*}
with $x=\|f(t)\|_{L^2(\Omega)},$ $y=\|u(t)\|_{L^2(\Omega)},$ $p=2\delta$ and $\eps=(\alpha\delta)^\frac1{2\delta},$ one obtains with Cauchy-Schwarz's inequality,
\begin{gather}
\label{demthmgenext4}
2\vint_\Omega|f(t,x)||u(t,x)|\d x
\le\frac{2\delta-1}\delta(\alpha\delta)^{-\frac1{2\delta-1}}\|f(t)\|_{L^2(\Omega)}^\frac{2\delta}{2\delta-1}
+\alpha\|u(t)\|_{L^2(\Omega)}^{2\delta}.
\end{gather}
Finally, set for any $t\ge0,$ $y(t)=\|u(t)\|_{L^2(\Omega)}^2$ and let us prove Property~\ref{pthmgenext1}). If $f$ satisfies \eqref{thmgenext2} then \eqref{demthmgenext3} may be rewritten as,
\begin{gather}
\label{odi}
y^\p(t)+2\alpha y(t)^\delta\le0,
\end{gather}
for almost every $t>T_0.$ We then conclude with the help of Lemma~\ref{lemodi}. Now assume that \eqref{thmgenext3}--\eqref{thmgenext4} hold where the constant $\eps_\star$ has to be determined later. We then have,
\begin{gather}
\label{demthmgenext5}
y(0)^{1-\delta}\le\alpha\,\delta\,(1-\delta)\,T_0,	\\
\label{demthmgenext6}
\|f(t)\|_{L^2(\Omega)}^2\le\eps_\star\|u\|_{L^\infty((0,\infty);H^\ell(\Omega))}^{-\frac{N(1-m)}{2\ell}\frac1{1-\delta}}\big(T_0-t\big)_+^\frac{2\delta-1}{1-\delta},
\end{gather}
where \eqref{demthmgenext5} is a consequence of \eqref{thmgenext3} and \eqref{demthmgenext6} is nothing else but \eqref{thmgenext4}. Gathering together \eqref{demthmgenext3}, \eqref{demthmgenext4} and \eqref{demthmgenext6}, one gets
\begin{gather*}
y^\p(t)+\alpha y(t)^\delta
\le\frac{2\delta-1}\delta(\Im(a)C_\GN^{-1}\delta)^{-\frac1{2\delta-1}}\eps_\star^\frac\delta{2\delta-1}
\|u\|_{L^\infty((0,\infty);H^\ell(\Omega))}^{-\frac{N(1-m)}{2\ell}\frac1{1-\delta}}\big(T_0-t\big)_+^\frac\delta{1-\delta}.
\end{gather*}
Choosing $\eps_\star=(2\delta-1)^{-\frac{2\delta-1}\delta}(\Im(a)C_\GN^{-1}\delta)^\frac1{1-\delta}(1-\delta)^\frac{2\delta-1}{\delta(1-\delta)},$ one obtains,
\begin{gather*}
y^\p(t)+\alpha y(t)^\delta\le y_\star\big(T_0-t\big)_+^\frac\delta{1-\delta}.
\end{gather*}
for almost every $t>0,$ where $y_\star$ is given by~\eqref{lemodif1}. Notice that \eqref{demthmgenext5} is nothing else but \eqref{lemodif3}. We infer by Lemma~\ref{lemodif} that $y(t)=0,$ for any $t\ge T_0.$
\medskip
\end{vproof}

\section{Proofs of the existence and uniqueness theorems}
\label{proofexi}

\begin{lem}
\label{Amaxmon}
Let Assumption~$\ref{ass}$ be fulfilled. Let us define the following $($nonlinear$)$ operator on $L^2(\Omega).$
\begin{gather}
\label{DA}
\begin{cases}
D(A)=\left\{u\in H^1_0(\Omega); \; \Delta u\in L^2(\Omega)\right\},	\medskip \\
\forall u\in D(A), \;  Au=-\vi\Delta u-\vi a|u|^{-(1-m)}u,
\end{cases}
\end{gather}
Then $A$ is a maximal monotone operator on $L^2(\Omega)$ $($and so $m$-accretive$)$ with dense domain.
\end{lem}

\noindent
The proof relies on the following lemmas.

\begin{lem}[\textbf{\cite{MR1224619}}]
\label{lemmon}
Let $0<m\le1.$ Set for any $z\in\C,$ $g(z)=|z|^{-(1-m)}z$ $(g(0)=0).$ Then for any $(z_1,z_2)\in\C\times\C,$
\begin{gather}
\label{lemmon1}
2\sqrt m\left|\Im\Big(\big(g(z_1)-g(z_2)\big)\big(\ovl{z_1-z_2}\big)\Big)\right|
		\le(1-m)\Re\Big(\big(g(z_1)-g(z_2)\big)\big(\ovl{z_1-z_2}\big)\Big),	\\
\label{m}
|g(z_1)-g(z_2)|\le3|z_1-z_2|^m.
\end{gather}
Let $\Omega\subseteq\R^N$ be an open subset. We define the mapping for any measurable function $u:\Omega\tends\C,$ which we still denote by $g,$ by $g(u)(x)=g(u(x)).$ Then for any $p\in[1,\infty),$
\begin{align}
\label{lemmon2}
&	g\in C\big(L^p(\Omega);L^\frac{p}m(\Omega)\big) \text{ and } g \text{ is bounded on bounded sets,}				\\
\label{lemmon3}
&	g\in C\big(L^2(\Omega);L^2(\Omega)\big) \text{ and } g \text{ is bounded on bounded sets, if } |\Omega|<\infty.
\end{align}
Finally, let $a\in\C$ with $\Im(a)>0$ satisfying~\eqref{a}. If $\big(g(u)-g(v)\big)(\ovl{u-v})\in L^1(\Omega)$ then,
\begin{gather}
\label{lemmon4}
\Re\left(-\vi\,a\vint_\Omega(g(u)-g(v))(\ovl{u-v})\d x\right)\ge0.
\end{gather}
We may choose, for instance, $u,v\in L^2(\Omega),$ if $|\Omega|<\infty,$ or $u,v\in L^{m+1}(\Omega),$ in the general case.
\end{lem}

\begin{proof*}
Estimate \eqref{lemmon1} is Lemma~2.2 of Liskevich and Perel$^\p$muter~\cite{MR1224619} while \eqref{m} comes from Lemma~\ref{lemm}, implying \eqref{lemmon2} and \eqref{lemmon3}. Finally, by \eqref{lemmon2}, \eqref{lemmon3} and Hölder's inequality, we have $\big(g(u)-g(v)\big)(\ovl{u-v})\in L^1(\Omega),$ for any $u,v$ as in the statement of the lemma and by \eqref{lemmon1},
\begin{align*}
	& \; \Re\left(-\vi\,a\vint_\Omega\big(g(u)-g(v)\big)(\ovl{u-v})\d x\right)								\\
   =	& \; \Im(a)\Re\vint_\Omega\big(g(u)-g(v)\big)\big(\ovl{u-v}\big)\d x
   		+\Re(a)\Im\vint_\Omega\big(g(u)-g(v)\big)\big(\ovl{u-v}\big)\d x								\\
 \ge	& \; \left(\Im(a)-|\Re(a)|\frac{1-m}{2\sqrt m}\right)\Re\vint_\Omega\big(g(u)-g(v)\big)\big(\ovl{u-v}\big)\d x	\\
 \ge	& \; 0.
\end{align*}
This ends the proof.
\medskip
\end{proof*}

\begin{vproof}{of Lemma~\ref{Amaxmon}.}
The density of the domain of the operator is obvious. Let $g$ be as in Lemma~\ref{lemmon}. It is well known that $(-\vi\Delta,D(A))$ is a maximal monotone operator on $L^2(\Omega)$ (Proposition~2.6.12, p.31, in Cazenave and Haraux~\cite{MR2000e:35003}). In addition, if we define $B$ on $L^2(\Omega)$ by $Bu=-\vi ag(u),$ it follows from \eqref{lemmon2}--\eqref{lemmon4} that $B\in C(L^2(\Omega);L^2(\Omega))$ and
\begin{gather*}
(Bu-Bv,u-v)_{L^2(\Omega)}=\Re\left(-\vi\,a\vint_\Omega(g(u)-g(v))(\ovl{u-v})\d x\right)\ge0,
\end{gather*}
for any $u,v\in L^2(\Omega).$ We then infer that $A=-\vi\Delta+B$ is a maximal monotone operator (Brezis~\cite{MR0348562}, Corollary~2.5, p.33 and Corollary~2.7, p.36).
\medskip
\end{vproof}

\noindent
To obtain \eqref{strongH12}, we need to regularize the nonlinearity in order to apply the $\nabla$ operator. We then establish the next lemma.

\begin{lem}
\label{lemenereg}
Let $\Omega\subseteq\R^N$ be an open subset, let $0<m<1,$ let $a\in\C$ with $\Im(a)>0$ satisfying \eqref{a} and let $\eps\in(0,1).$ Let for any $u\in L^2(\Omega),$ $g_\eps(u)=(|u|^2+\eps)^{-\frac{1-m}2}u.$ Finally, let $g$ be as in Lemma~$\ref{lemmon}$ and let $D(A)$ be defined by \eqref{DA}. Then,
\begin{align}
\label{lemenereg1}
&	g_\eps\in C\big(L^2(\Omega);L^2(\Omega)\big)\cap C\big(H^1_0(\Omega);H^1_0(\Omega)\big),	\\
\label{lemenereg2}
&	\forall u\in D(A), \; \Re\left(\vi a\vint_\Omega g_\eps(u)\ovl{\Delta u}\d x\right)\ge0,				\\
\label{lemenereg3}
&	\forall u\in D(A) \text{ such that } u^m\Delta u\in L^1(\Omega), \; \Re\left(\vi a\vint_\Omega g(u)\ovl{\Delta u}\d x\right)\ge0.
\end{align}
\end{lem}

\begin{rmk}
\label{rmklemenereg}
If $\Omega\subseteq\R^N$ is arbitrary, $m=1$ and $\Im(a)>0$ then for any $u\in D(A),$
\begin{gather*}
\Re\left(\vi a\int_\Omega g(u)\ovl{\Delta u}\d x\right)=\Im(a)\|\nabla u\|_{L^2(\Omega)}^2\ge0.
\end{gather*}
In other words, one directly obtains \eqref{lemenereg3}.
\end{rmk}

\begin{vproof}{of Lemma~\ref{lemenereg}.}
A straightforward calculation shows that for any $\eps\in(0,1),$
\begin{gather*}
|g_\eps(u)-g_\eps(v)|\le C\eps^{-1}|u-v|,	\\
|\nabla g_\eps(u)|\le C\eps^{-1}|\nabla u|.
\end{gather*}
It follows that if $u\in H^1_0(\Omega)$ then $g_\eps(u)\in H^1_0(\Omega)$ and \eqref{lemenereg1} comes from the above estimates and the partial converse of the dominated convergence Theorem (see, for instance, Brezis~\cite{MR2759829}, Theorem~4.9, p.94). Let us turn out to the proof of  \eqref{lemenereg2}. Let $u\in D(A).$ It follows from \eqref{lemenereg1} that we can take the scalar product in $L^2$ between $\vi ag_\eps(u)$ and $\Delta u.$ We then obtain,
\begin{align*}
	&	\; \Re\left(\vi a\vint_\Omega g_\eps(u)\ovl{\Delta u}\d x\right)=
			(\vi ag_\eps(u),\Delta u)_{L^2(\Omega)}=-(\vi a\nabla g_\eps(u),\nabla u)_{L^2(\Omega)}				\\
   =	&	\; \Re\left(-\vi a\vint_\Omega\frac{|\nabla u|^2(|u|^2+\eps)-(1-m)\Re(u\ovl{\nabla u}).u\ovl{\nabla u}}
			{(|u|^2+\eps)^\frac{3-m}2}\d x\right)														\\
   =	&	\; \Im(a)\vint_\Omega\frac{|\nabla u|^2(|u|^2+\eps)-(1-m)|\Re(u\ovl{\nabla u})|^2}{(|u|^2+\eps)^\frac{3-m}2}\d x
			-\Re(a)\vint_\Omega\frac{(1-m)\Re(u\ovl{\nabla u}).\Im(u\ovl{\nabla u)}}{(|u|^2+\eps)^\frac{3-m}2}\d x		\\
   =	&	\; \eps\,\Im(a)\vint_\Omega\frac{|\nabla u|^2}{(|u|^2+\eps)^\frac{3-m}2}\d x								\\
	&	\; +\Im(a)\vint_\Omega\frac{m|\Re(u\ovl{\nabla u})|^2
			+|\Im(u\ovl{\nabla u})|^2}{(|u|^2+\eps)^\frac{3-m}2}\d x
			-\Re(a)\vint_\Omega\frac{(1-m)\Re(u\ovl{\nabla u}).\Im(u\ovl{\nabla u)}}{(|u|^2+\eps)^\frac{3-m}2}\d x,
\end{align*}
where we used in the last equality the fact that, $|\nabla u|^2|u|^2=|\Re(u\ovl{\nabla u})|^2+|\Im(u\ovl{\nabla u})|^2.$ To conclude, it remains to show that,
\begin{gather}
\label{demlemenereg}
(1-m)|\Re(a)|\,|\Re(u\ovl{\nabla u})|\,|\Im(u\ovl{\nabla u})|\le \Im(a)\left(m|\Re(u\ovl{\nabla u})|^2+|\Im(u\ovl{\nabla u})|^2\right).
\end{gather}
Using our assumption on $a$ and the following Young inequality,
\begin{gather*}
2|xy|\le\delta x^2+\frac{y^2}\delta,
\end{gather*}
with $x=|\Re(u\ovl{\nabla u})|,$ $y=|\Im(u\ovl{\nabla u})|$ and $\delta=\sqrt m,$ we obtain,
\begin{align*}
	&	\; (1-m)|\Re(a)|\,|\Re(u\ovl{\nabla u})|\,|\Im(u\ovl{\nabla u)}|								\\
  \le	&	\; 2\sqrt m\,\Im(a)|\Re(u\ovl{\nabla u})|\,|\Im(u\ovl{\nabla u})|							\\
  \le	&	\; \sqrt m\,\Im(a)\left(\sqrt m|\Re(u\ovl{\nabla u})|^2+\frac{|\Im(u\ovl{\nabla u})|^2}{\sqrt m}\right)	\\
  \le	&	\;  \Im(a)\left(m|\Re(u\ovl{\nabla u})|^2+|\Im(u\ovl{\nabla u})|^2\right),
\end{align*}
which is \eqref{demlemenereg}. Finally, since we have $g_\eps(u)\xrightarrow[\eps\searrow0]{\text{a.e. on }\Omega}g(u)$ and $|g_\eps(u)|\overset{\text{a.e.}}{\le}|g(u)|,$ for any $\eps>0,$ \eqref{lemenereg3} is a consequence of \eqref{lemenereg2} and the dominated convergence Theorem.
\medskip
\end{vproof}

\noindent
Concerning the continuous dependence with respect to the data we have:

\begin{lem}
\label{lemdep}
Let $\Omega\subseteq\R^N$ be an open subset, $0<m\le1$ and $a\in\C$ with $\Im(a)>0$ satisfying~\eqref{a}. Let $X=L^2(\Omega)\cap L^{m+1}(\Omega)$ or $X=H^1_0(\Omega)\cap L^{m+1}(\Omega).$ Finally, let $f_1,f_2\in L^1_\loc([0,\infty);L^2(\Omega))$ and let
\begin{gather*}
u,v\in L^p_\loc\big([0,\infty);X\big)\cap W^{1,p^\p}_\loc\big([0,\infty);X^\star\big),
\end{gather*}
 for some $1<p<\infty.$ If,
\begin{gather*}
\vi u_t+\Delta u+a|u|^{-(1-m)}u=f_1,	\\
\vi v_t+\Delta v+a|v|^{-(1-m)}v=f_2,
\end{gather*}
in $\Dr^\p\big((0,\infty)\times\Omega\big),$ then $u,v\in C\big([0,\infty);L^2(\Omega)\big)$ and
\begin{gather}
\label{lemdep1}
\|u(t)-v(t)\|_{L^2(\Omega)}\le\|u(s)-v(s)\|_{L^2(\Omega)}+\vint_s^t\|f_1(\sigma)-f_2(\sigma)\|_{L^2(\Omega)}\d\sigma,
\end{gather}
for any $t\ge s\ge0.$
\end{lem}

\begin{proof*}
By Lemma~\ref{lemden0} and the dense embedding $X\inj L^2(\Omega),$ we have $L^2(\Omega)\inj X^\star\inj\Dr^\p\Omega)$ and for any $(x,y)\in L^2(\Omega)\times X,$
\begin{gather}
\label{demlemdep}
(x,y)_{L^2(\Omega)}=\langle x,y\rangle_{L^2(\Omega),L^2(\Omega)}=\langle x,y\rangle_{X^\star,X}.
\end{gather}
It follows from above and \eqref{gmm} that the equations in the lemma make sense in $X^\star$ and we then have,
\begin{gather*}
\vi(u-v)_t+\Delta(u-v)+\big(ag(u)-ag(v)\big)=f_1-f_2, \; \text{ in } X^\star,
\end{gather*}
almost everywhere on $(0,\infty),$ where $g$ is as in Lemma~\ref{lemmon}. Taking the $X^\star-X$ duality product of the above equation with $\vi(u-v),$ it follows from \ref{embL21}) of Lemma~\ref{lemC}, \ref{lemmass}) of Lemma~\ref{lemmassene} and \eqref{demlemdep} that $u,v\in C\big([0,\infty);L^2(\Omega)\big),$ the mapping $t\longmapsto\|u(t)-v(t)\|_{L^2(\Omega)}^2$ belongs to $W^{1,1}_\loc\big([0,\infty);\R\big)$ and,
\begin{gather*}
\frac12\frac\d{\d t}\|u(\:.\:)-v(\:.\:)\|_{L^2(\Omega)}^2
+\big\langle ag(u)-ag(v),\vi(u-v)\big\rangle_{X^\star,X}=\big(f_1-f_2,\vi(u-v)\big)_{L^2(\Omega)},
\end{gather*}
almost everywhere on $(0,\infty).$ Applying \eqref{dualg}, \eqref{lemmon4} and Cauchy-Schwarz's inequality to the above, one infers
\begin{gather*}
\frac12\frac\d{\d t}\|u(\:.\:)-v(\:.\:)\|_{L^2(\Omega)}^2\le\|f_1-f_2\|_{L^2(\Omega)}\|u-v\|_{L^2(\Omega)},
\end{gather*}
almost everywhere on $(0,\infty).$ Integrating over $(s,t),$ one obtains \eqref{lemdep1}.
\medskip
\end{proof*}

\begin{vproof}{of Theorem~\ref{thmstrongH2}.}
By  Lemma~\ref{Amaxmon} and Vrabie~\cite{MR1375237} (Theorem~1.7.1, p.23), there exists a unique $u\in W^{1,\infty}_\loc\big([0,\infty);L^2(\Omega)\big)$ satisfying $u(t)\in H^1_0(\Omega),$ $\Delta u(t)\in L^2(\Omega)$ and \eqref{nls} in $L^2(\Omega),$ for almost every $t>0,$ $u(0)=u_0$ and \eqref{strongH23}. Then \eqref{strongH21} comes from \eqref{strongH23}. It follows from \ref{W11}) of Lemma~\ref{lemC}, \eqref{lemmon2}--\eqref{lemmon3}, \eqref{strongH23}, \eqref{E} and \eqref{nls} that,
\begin{gather}
\label{demH21}
f\in C\big([0,\infty);L^2(\Omega)\big),					\\
\label{demH22}
|u|^{-(1-m)}u\in C\big([0,\infty);L^2(\Omega)\big),		\\
\label{demH23}
\Delta u\in L^\infty_\loc\big([0,\infty);L^2(\Omega)\big),	\\
\nonumber
u\in L^\infty_\loc\big([0,\infty);H^1_0(\Omega)\big),
\end{gather}
so that $u$ is a $H^2$-solution and $u\in C\big([0,\infty);H^1_0(\Omega)\big)$ (by \ref{embL22}) of Lemma~\ref{lemC}). So,
\begin{gather}
\label{demH24}
\Delta u\in C\big([0,\infty);H^{-1}(\Omega)\big).
\end{gather}
It then follows from \eqref{demH21}, \eqref{demH22}, \eqref{demH24} and \eqref{nls} that,
\begin{gather*}
u_t\in C\big([0,\infty);H^{-1}(\Omega)\big).
\end{gather*}
By \eqref{E}, \eqref{strongH21} and \eqref{demH23}, one obtains \eqref{strongH22} and Properties~\ref{thmstrongH21}) and \ref{thmstrongH22}) are proved. Property~\ref{thmstrongH23}) follows easily from Property~\ref{thmstrongH21}), \eqref{M'} and Proposition~\ref{propL2}. Finally, Property~\ref{thmstrongH24}) comes from \eqref{lemdep1}, \eqref{strongH23}, \eqref{E}, \eqref{lemmon2}, \eqref{lemmon3}, the embedding \ref{W11}) of Lemma~\ref{lemC} and \eqref{nls}. This concludes the proof of the theorem.
\medskip
\end{vproof}

\begin{vproof}{of Theorem~\ref{thmweak}.}
Existence comes from density of $H^2_0(\Omega)\times W^{1,1}_\loc([0,\infty);L^2(\Omega))$ in $L^2(\Omega)\times L^1_\loc([0,\infty);L^2(\Omega)),$ Theorem~\ref{thmstrongH2}, \eqref{lemdep1} and completeness of $C\big([0,T];L^2(\Omega)\big),$ for any $T>0.$ Property~\ref{thmweak1}) comes from Proposition~\ref{propL2}. Estimate~\eqref{estthmweak} being stable by passing to the limit in $C\big([0,T];L^2(\Omega)\big)\times L^1\big((0,T);L^2(\Omega)\big),$ for any $T>0,$ it is sufficient to establish it for the $H^2$-solutions. This then comes from Lemma~\ref{lemdep} and the uniqueness conclusion of the theorem follows. Finally, Property~\ref{thmweak1}) comes from Proposition~\ref{propL2}.
\medskip
\end{vproof}

\begin{vproof}{of Theorem~\ref{thmstrongH1}.}
The uniqueness of solutions comes from Lemma~\ref{lemdep}. Let $f\in W^{1,1}_\loc([0,\infty);H^1_0(\Omega))$ and let $u_0\in H^1_0(\Omega).$ Let $(\vphi_n)_{n\in\N}\subset H^2_0(\Omega)$ be such that $\vphi_n\xrightarrow[n\to\infty]{H^1_0(\Omega)}u_0.$ Finally, let $g$ be defined as in Lemma~\ref{lemmon} and for each $n\in\N,$ let $u_n$ be the unique $H^2$-solution of \eqref{nls}--\eqref{nlsb} such that $u_n(0)=\vphi_n,$ given by Theorem~\ref{thmstrongH2}. By Lemma~\ref{lemdep}, we have for any $T>0$ and $n,p\in\N,$
\begin{gather}
\label{demthmstrongH11}
\|u_n\|_{C([0,T];L^2(\Omega))}\le\|\vphi_n\|_{L^2(\Omega)}+\int_0^T\|f(t)\|_{L^2(\Omega)}\d t,	\\
\nonumber
\|u_n-u_p\|_{L^\infty((0,\infty);L^2(\Omega))}\le\|\vphi_n-\vphi_p\|_{L^2(\Omega)},
\end{gather}
It follows that for any $T>0,$ $(u_n)_{n\in\N}$ is a Cauchy sequence in $C\big([0,T];L^2(\Omega)\big).$ As a consequence, and with \eqref{lemmon2}--\eqref{lemmon3}, there exists $u\in C\big([0,\infty);L^2(\Omega)\big)$ such that for any $T>0,$
\begin{gather}
\label{demthmstrongH12}
u_n\xrightarrow[n\to\infty]{C([0,T];L^2(\Omega))}u,		\\
\label{demthmstrongH13}
g(u)\in C\big([0,T];L^2(\Omega)\big),					\\
\label{demthmstrongH14}
g(u_n)\xrightarrow[n\to\infty]{C([0,T];L^2(\Omega))}g(u).
\end{gather}
By definition, it follows from \eqref{demthmstrongH12} that $u$ is a weak solution of \eqref{nls}--\eqref{u0} (take $f_n=f,$ for any $n\in \N).$ By \ref{rmkdefsol3}) of Remark~\ref{rmkdefsol}, we can take the $L^2$-scalar product of \eqref{nls} with $-\vi\Delta u_n$ and it follows from \eqref{E'} that for any $n\in\N$ and almost every $s>0,$
\begin{gather*}
\frac12\frac{\d}{\d t}\|\nabla u_n(s)\|_{L^2(\Omega)}^2+\Re\left(\vi a\vint_\Omega g(u_n(s))\ovl{\Delta u_n(s)}\d x\right)
=\big(\nabla f(s),\vi\nabla u_n(s)\big)_{L^2(\Omega)},
\end{gather*}
which gives with \eqref{lemenereg3}, Remark~\ref{rmklemenereg} and Cauchy-Schwarz's inequality,
\begin{gather*}
\frac12\frac{\d}{\d t}\|\nabla u_n(s)\|_{L^2(\Omega)}^2\le\|\nabla f_n(s)\|_{L^2(\Omega)}\|\nabla u_n(s)\|_{L^2(\Omega)}.
\end{gather*}
By integration, we obtain for any $t>0$ and any $n\in\N,$ 
\begin{gather}
\label{demthmstrongH15}
\|\nabla u_n(t)\|_{L^2(\Omega)}\le\|\nabla \vphi_n\|_{L^2(\Omega)}+\int_0^t\|\nabla f(s)\|_{L^2(\Omega)}\d s.
\end{gather}
By the Sobolev embedding \ref{W11}) of Lemma~\ref{lemC},
\begin{gather}
\label{demthmstrongH16}
W^{1,1}_\loc\big([0,\infty);L^2(\Omega)\big)\inj C\big([0,\infty);L^2(\Omega)\big),
\end{gather}
\eqref{demthmstrongH11}, \eqref{demthmstrongH14}, \eqref{demthmstrongH15} and \eqref{nls}, we infer that,
\begin{gather}
\label{demthmstrongH17}
(u_n)_{n\in\N} \text{ is bounded in } L^\infty\big((0,T);H^1_0(\Omega)\big)\cap W^{1,\infty}\big((0,T);H^{-1}(\Omega)\big),
\end{gather}
for any $T>0.$ Applying Proposition~1.3.14, p.12, and Proposition~1.1.2, p.2, in Cazenave~\cite{MR2002047}, it follows from \eqref{demthmstrongH12} and \eqref{demthmstrongH17} that,
\begin{align}
\label{demthmstrongH18}
&	u\in C_\w\big([0,\infty);H^1_0(\Omega)\big)\cap W^{1,\infty}_\loc\big([0,\infty);H^{-1}(\Omega)\big), 	\\
\label{demthmstrongH19}
&	\Delta u\in C\left([0,\infty);H^{-2}(\Omega)\right),											\\
\label{demthmstrongH110}
&	u_n(t)\weak u(t), \; \text{ in } \; H^1_\w(\Omega), \; \text{ as } \; n\to\infty,
\end{align}
for any $t\ge0.$ Since $u$ is a weak solution, $u$ solves \eqref{nls} in $H^{-2}(\Omega),$ for almost every $t>0$ (Property~\ref{rmkdefsol5}) of Remark~\ref{rmkdefsol}). As a consequence, and with help of \eqref{demthmstrongH13}, \eqref{demthmstrongH16} and \eqref{demthmstrongH19}, we have that $u_t\in C\big([0,\infty);H^{-2}(\Omega)\big)$ and $u$ satisfies \eqref{nls} in $H^{-2}(\Omega),$ for any $t\ge0.$ We then infer with \eqref{demthmstrongH18} that $u$ is a $H^1_0$-solution and Property~\ref{thmstrongH11}) holds. Still by \eqref{demthmstrongH18}, we have for any $t\ge s\ge0,$
\begin{align*}
	&	\; \|u(t)-u(s)\|_{L^2(\Omega)}^2\le2\|u\|_{L^\infty((s,t);H^1_0(\Omega))}\|u(t)-u(s)\|_{H^{-1}(\Omega)}	\\
  \le	&	\; 2\|u\|_{L^\infty((s,t);H^1_0(\Omega))}\|u_t\|_{L^\infty((s,t);H^{-1}(\Omega))}|t-s|,
\end{align*}
which is \eqref{strongH11}. By \eqref{demthmstrongH110}, the weak lower semicontinuity of the norm and \eqref{demthmstrongH15}, one obtains \eqref{strongH12} and Property~\ref{thmstrongH12}) is proved. Property~\ref{thmstrongH13}) follows easily from Proposition~\ref{propL2} and the fact that $u,f\in C\big([0,\infty);L^2(\Omega)\big)$ and $L^2(\Omega)\inj L^{m+1}(\Omega).$ Finally, Property~\ref{thmstrongH14}) comes from \eqref{estthmweak}, \eqref{strongH12}, \eqref{lemmon2}, \eqref{lemmon3}, \ref{W11}) of Lemma~\ref{lemC} and \eqref{nls}. This concludes the proof of the theorem.
\medskip
\end{vproof}

\section{Proofs of the finite time extinction property and asymptotic behavior theorems}
\label{proofext}

\begin{vproof}{of Theorem~\ref{thmextH2}.}
For the Property~$\ref{thmextH21}),$ apply Theorems~\ref{thmstrongH1}, \ref{thmstrongH2}, Remark~\ref{rmkthmstrong} and Theorem~\ref{thmgenext} (with $\ell=1,$ if $u_0\in H^1_0(\Omega)$ and $\ell=2,$ if $u_0\in H^2(\Omega)\cap H^1_0(\Omega)).$ We then obtain the finite time extinction result and the upper bound on $T_\star.$ The lower bound on $T_\star$ comes from \ref{rmkthmgenext2}) of Remark~\ref{rmkthmgenext}. Property~$\ref{thmextH22})$ comes from Remark~\ref{rmkthmstrongH10}.
\medskip
\end{vproof}

\begin{vproof}{of Theorem~\ref{thm0s}.}
By Theorems~\ref{thmstrongH1}, \ref{thmstrongH2} and Remark~\ref{rmkthmstrong}, $u\in L^\infty\big((0,\infty);H^\ell(\Omega)\big),$ where $\ell=1,$ if $u_0\in H^1_0(\Omega)$ and $\ell=2,$ if $u_0\in H^2(\Omega)\cap H^1_0(\Omega).$ The result then comes from \ref{rmkthmgenext3}) of Remark~\ref{rmkthmgenext}.
\medskip
\end{vproof}

\begin{vproof}{of Theorem~\ref{thm0w}.}
Let the assumptions of the theorem be fulfilled. We proceed to the proof in two steps. \\
\textbf{Step 1.} Assume further that $f\in\Dr\big([0,\infty);L^2(\Omega)\big)$ and $u_0\in H^2_0(\Omega).$ Then, $\vlim_{t\nearrow\infty}\|u(t)\|_{L^2(\Omega)}=0.$ \\
It follows from uniqueness and Theorem~\ref{thmstrongH2} that $u$ is a $H^2$-solution and $u\in L^\infty\big((0,\infty);H^1_0(\Omega)\big).$ Let $[0,T_0]\supset\supp f.$ By \eqref{L2}, $\frac{\d}{\d t}\|u(t)\|_{L^2(\Omega)}^2\le0,$ for any $t>T_0.$ It follows that $\vlim_{t\nearrow\infty}\|u(t)\|_{L^2(\Omega)}=\ell_0,$ for some $\ell_0\in[0,\infty).$ If $m=1$ then we have, one more time by \eqref{L2}, $\frac{\d}{\d t}\|u(t)\|_{L^2(\Omega)}^2\le-2\Im(a)\ell_0^2,$ for any $t>T_0.$ It follows that $\ell_0=0.$ Now, assume that $m<1$ and suppose, by contradiction, that $\ell_0\neq0.$ Let $q\in(2,\infty)$ with $(N-2)q<2N.$ By Hölder's inequality and Sobolev's embedding $H^1_0(\Omega)\inj L^q(\Omega),$ there exists $\theta\in(0,1)$ such that,
\begin{gather*}
0<\ell_0\le\|u(t)\|_{L^2(\Omega)}\le\|u(t)\|_{L^{m+1}(\Omega)}^\theta\|u(t)\|_{L^q(\Omega)}^{1-\theta}\le C\|u(t)\|_{L^{m+1}(\Omega)}^\theta\|u\|_{L^\infty((0,\infty);H^1_0(\Omega))}^{1-\theta},
\end{gather*}
for any $t>T_0.$ We infer that, $\vinf_{t>T_0}\|u(t)\|_{L^{m+1}(\Omega)}>0,$ which implies with \eqref{L2},
\begin{gather*}
\frac{\d}{\d t}\|u(t)\|_{L^2(\Omega)}^2\le-2\Im(a)\vinf_{t>T_0}\|u(t)\|_{L^{m+1}(\Omega)}^{m+1}<0,
\end{gather*}
for any $t>T_0.$ As a consequence, $\vlim_{t\nearrow\infty}\|u(t)\|_{L^2(\Omega)}=-\infty,$ a contradiction.
\\
\textbf{Step 2.} Conclusion. \\
Let $(\vphi_n)_{n\in\N}\subset H^2_0(\Omega)$ and $(f_n)_{n\in\N}\subset\Dr\big([0,\infty);L^2(\Omega)\big)$ be such that, 
\begin{gather*}
\vphi_n\xrightarrow[n\to\infty]{L^2(\Omega)}u_0 \; \text{ and } \;
f_n\xrightarrow[n\to\infty]{L^1((0,\infty);L^2(\Omega))}f.
\end{gather*}
For each $n\in\N,$ let $u_n$ the $H^2$-solution to \eqref{nls}--\eqref{nlsb}, with $f_n$ instead of $f,$ be such that $u_n(0)=\vphi_n,$ given by Theorem~\ref{thmstrongH2}. Let $n\in\N.$ It follows from \eqref{estthmweak} that,
\begin{align*}
	&	\; \|u(t)\|_{L^2(\Omega)}\le\|u-u_n\|_{L^\infty((0,\infty);L^2(\Omega))}+\|u_n(t)\|_{L^2(\Omega)}		\\
  \le	&	\; \|u_0-\vphi_n\|_{L^2(\Omega)}+\|f-f_n\|_{L^1((0,\infty);L^2(\Omega))}+\|u_n(t)\|_{L^2(\Omega)},
\end{align*}
for any $t>0.$ We get from Step~1,
\begin{gather*}
\limsup_{t\nearrow\infty}\|u(t)\|_{L^2(\Omega)}\le\|u_0-\vphi_n\|_{L^2(\Omega)}+\|f-f_n\|_{L^1((0,\infty);L^2(\Omega))}.
\end{gather*}
Letting $n\nearrow\infty,$ we obtain $\vlim_{t\nearrow\infty}\|u(t)\|_{L^2(\Omega)}=0.$ Finally, the general case comes from the embedding $L^2(\Omega)\inj L^p(\Omega),$ which holds for any $p\in(0,2],$ as soon as $|\Omega|<\infty.$ This concludes the proof.
\medskip
\end{vproof}

\appendix
\section{Appendix}
\label{appendix}

In this appendix, we recall some useful estimates and results about Sobolev spaces.

\begin{lem}
\label{lemm}
Let $0<m\le1.$ Then we have for any $(z_1,z_2)\in\C\times\C,$
\begin{gather}
\left||z_1|^{-(1-m)}z_1-|z_2|^{-(1-m)}z_2\right|\le3|z_1-z_2|^m,
\end{gather}
where $|z|^{-(1-m)}z=0,$ if $z=0.$
\end{lem}

\begin{proof*}
Let $0<m<1$ (the case $m=1$ being obvious). We proceed to the proof in four steps. \\
\textbf{Step 1:} $\forall t,s\ge0,$ $|t^m-s^m|\le|t-s|^m.$ \\
Let for $x\ge1,$ $f(x)=(x-1)^m-(x^m-1).$ Then $f^\p>0$ on $(1,\infty)$ and so $f\left(\frac{t}s\right)\ge f(1)=0,$ for any $t\ge s>0.$ Hence Step~1. \\
\textbf{Step 2:} $\forall a\ge0,$ $\forall\theta\in\R,$ $\left|a^m-a^me^{\vi\theta}\right|\le2^{1-m}\left|a-ae^{\vi\theta}\right|^m.$ \\
We have for any $\theta\in\R,$ $\left|1-e^{\vi\theta}\right|^{1-m}\le2^{1-m},$ implying $\left|1-e^{\vi\theta}\right|\le2^{1-m}\left|1-e^{\vi\theta}\right|^m,$ therefore Step~2. \\
\textbf{Step 3:} $\forall(z_1,z_2)\in\C\setminus\{0\}\times\C,$ $\left||z_2|-\frac{\ovl{z_1}}{|z_1|}z_2\right|^m\le2^m|z_1-z_2|^m.$ \\
We have,
\begin{align*}
	&	\; \left||z_2|-\frac{\ovl{z_1}}{|z_1|}z_2\right|
			=\left|\left(|z_2|-\frac{\ovl{z_1}}{|z_1|}z_1\right)+\left(\frac{\ovl{z_1}}{|z_1|}z_1
			-\frac{\ovl{z_1}}{|z_1|}z_2\right)\right|	\\
   =	&	\; \left|\big(|z_2|-|z_1|\big)+\left(\frac{\ovl{z_1}}{|z_1|}z_1-\frac{\ovl{z_1}}{|z_1|}z_2\right)\right|
			\le\big||z_2|-|z_1|\big|+|z_1-z_2|\le2|z_1-z_2|.
\end{align*}
Hence Step~3. \\
\textbf{Step 4:} Conclusion. \\
Let $(z_1,z_2)\in\C\times\C$ with $z_1z_2\neq0,$ otherwise there is nothing to prove.
\begin{align*}
	&	\; \left||z_1|^{-(1-m)}z_1-|z_2|^{-(1-m)}z_2\right|
			=\left||z_1|^{-(1-m)}z_1\frac{\ovl{z_1}}{|z_1|}-|z_2|^{-(1-m)}z_2\frac{\ovl{z_1}}{|z_1|}\right|			\\
   =	&	\; \left|\big(|z_1|^m-|z_2|^m\big)+\left(|z_2|^m-|z_2|^m\frac{\ovl{z_1}}{|z_1|}\frac{z_2}{|z_2|}\right)\right|
			\overset{\text{Steps~1 and 2}}{\le}
			|z_1-z_2|^m+2^{1-m}\left||z_2|-|z_2|\frac{\ovl{z_1}}{|z_1|}\frac{z_2}{|z_2|}\right|^m				\\
   =	&	\; |z_1-z_2|^m+2^{1-m}\left||z_2|-\frac{\ovl{z_1}}{|z_1|}z_2\right|^m
			\overset{\text{Steps~3}}{\le}3|z_1-z_2|^m.
\end{align*}
The lemma is proved.
\medskip
\end{proof*}

\noindent
The four next lemmas are, more or less, a repetition of some similar results contained in the unpublished book by Brezis and Cazenave~\cite{bc}.

\begin{lem}
\label{lemden0}
Let $\Omega\subseteq\R^N$ be a nonempty open subset, let $k,m\in\N_0$ and let $1\le p,q<\infty.$ Then $\Dr(\Omega)\inj W^{k,p}_0(\Omega)\cap W^{m,q}_0(\Omega)$ with dense embedding. In addition, $W^{k,p}_0(\Omega)\cap W^{m,q}_0(\Omega)$ is separable and,
\begin{gather}
\label{dual}
\big(W^{k,p}_0(\Omega)\cap W^{m,q}_0(\Omega)\big)^\star=W^{-k,p^\p}(\Omega)+W^{-m,q^\p}(\Omega)\inj\Dr^\p(\Omega).
\end{gather}
Finally, if $p,q>1$ then $W^{k,p}_0(\Omega)\cap W^{m,q}_0(\Omega)$ and $W^{-k,p^\p}(\Omega)+W^{-m,q^\p}(\Omega)$ are reflexive and separable.
\end{lem}

\begin{proof*}
Set $X=W^{k,p}_0(\Omega)\cap W^{m,q}_0(\Omega).$ Without loss of generality, we may assume that $p\le q.$ It is clear that $\Dr(\Omega)\inj X.$ The equality in \eqref{dual} comes from the density of $\Dr(\Omega)$ in the spaces $W^{j,r}_0(\Omega)$ and Bergh and Löfström~\cite{MR0482275} (Lemma~2.3.1, p.24-25, and Theorem~2.7.1, p.32). Since for any $j\in\N_0$ and $r\in[1,\infty),$ $W^{-j,r^\p}(\Omega)\inj\Dr^\p(\Omega),$ we have by the equality in \eqref{dual},
\begin{gather*}
X^\star=\left\{T\in\Dr^\p(\Omega);T=T_1+T_2, (T_1,T_2)\in W^{-k,p^\p}(\Omega)\times W^{-m,q^\p}(\Omega)\right\}.
\end{gather*}
Let $T\in X^\star$ be such that $\langle T,\vphi\rangle_{X^\star,X}=0,$ for any $\vphi\in\Dr(\Omega).$ It follows from above that for any $\vphi\in\Dr(\Omega),$ $\langle T,\vphi\rangle_{\Dr^\p(\Omega),\Dr(\Omega)}=\langle T,\vphi\rangle_{X^\star,X}=0.$ Then $T=0$ in $\Dr^\p(\Omega),$ hence in $X^\star.$ We deduce that $\Dr(\Omega)\inj X$ is dense (Brezis~\cite{MR2759829}, Corollary~1.8, p.8) and so $X^\star\inj\Dr^\p(\Omega).$ Now, let $n>k+m$ be large enough to have $W^{n,p}_0(\Omega)\inj X.$ Since this embedding is dense and $W^{n,p}_0(\Omega)$ is separable, we infer that $X$ is separable. Finally, separability and reflexivity of the last part of the lemma present no difficulty and follow easily from reflexivity and separability of the spaces $W^{j,r}_0(\Omega),$ \eqref{dual} and Eberlein–\v{S}mulian's Theorem (Brezis~\cite{MR2759829}, Theorem~3.19, p.70, and Corollary~3.27, p.73).
\medskip
\end{proof*}

\begin{lem}[\cite{bc}]
\label{lemden}
Let $I\subseteq\R$ be an open interval, let $1\le p,q<\infty$ and let $X\inj Y$ be two Banach spaces. Then $\Dr(\ovl I;X)$ is dense in $L^p(I;X)\cap W^{1,q}(I;Y).$ Moreover, if $Z$ is a Banach space such that $Z\inj X$ with dense embedding then $\Dr(\ovl I;Z)$ is dense in $L^p(I;X)\cap W^{1,q}(I;Y).$
\end{lem}

\begin{proof*}
We first construct a linear extension operator to bring back to the case $I=\R.$  The first statement then follows from the standard procedure of truncation and regularization, while the second statement comes from the density of $\Dr(\R;Z)$ in $C^1_\co(\R;X),$ for the norm of $C^1_\b(\R;X).$
\medskip
\end{proof*}

\begin{lem}
\label{lemC}
Let $\Omega\subseteq\R^N$ be an open subset. Consider the Hilbert space given by $D(A)$ with,
\begin{gather*}
D(A)=\left\{u\in H^1_0(\Omega); \; \Delta u\in L^2(\Omega)\right\},		\\
\|u\|_{D(A)}^2=\|u\|_{H^1_0(\Omega)}^2+\|\Delta u\|_{L^2(\Omega)}^2,
\end{gather*}
for any $u\in D(A).$ Moreover, let $X$ be a Banach space, let $I$ be an open interval and let $1<p<\infty.$ We have the following results.
\begin{enumerate}[$1)$]
\item
\label{W11}
$W^{1,1}\big(I;X\big)\inj C_{\b,\vu}\big(\ovl I;X\big).$
\item
\label{embL21}
$L^p(I;X)\cap W^{1,p^\p}(I;X^\star)\inj C_\b\big(\ovl I;L^2(\Omega)\big),$ if $X\inj L^2(\Omega)$ with dense embedding.
\item
\label{embL22}
$L^p\big(I;D(A)\big)\cap W^{1,p^\p}\big(I;L^2(\Omega)\big)\inj C_\b\big(\ovl I;H^1_0(\Omega)\big).$
\end{enumerate}
\end{lem}

\begin{lem}
\label{lemmassene}
Let $\Omega\subseteq\R^N$ be an open subset, let $I$ be an open interval and let $1<p<\infty.$ For $t\in I$ and $u=u(t,x)\in\C,$ let us define $($formally$),$
\begin{gather*}
M(t)=\dfrac12\|u(t)\|_{L^2(\Omega)}^2 \; \text{ and } \; E(t)=\dfrac12\|\nabla u(t)\|_{L^2(\Omega)}^2.
\end{gather*}
Let $D(A)$ the Hilbert space be defined in Lemma~$\ref{lemC}$ and let $X\inj L^2(\Omega)$ be a Banach space with dense embedding. We then have the following results.
\begin{enumerate}[$1)$]
\item
\label{lemmass}
If $u\in L^p(I;X)\cap W^{1,p^\p}(I;X^\star)$ or if $u\in W^{1,1}(I;L^2(\Omega))$ then $M\in W^{1,1}(I;\R)$ and,
\begin{gather}
\label{M'}
M^\p(t)=
\begin{cases}
\big\langle u(t),u^\p(t)\big\rangle_{X,X^\star},	&	\text{if } \;	u\in L^p(I;X)\cap W^{1,p^\p}(I;X^\star),	\medskip \\
\big(u(t),u^\p(t)\big)_{L^2(\Omega)},			&	\text{if } \;	u\in W^{1,1}(I;L^2(\Omega)),
\end{cases}
\end{gather}
for almost every $t\in I.$
\item
\label{lemene}
If $u\in L^p(I;D(A))\cap W^{1,p^\p}(I;L^2(\Omega))$ then $E\in W^{1,1}(I;\R)$ and,
\begin{gather}
\label{E'}
E^\p(t)=\big(-\Delta u(t),u^\p(t)\big)_{L^2(\Omega)},
\end{gather}
for almost every $t\in I.$
\end{enumerate}
\end{lem}

\begin{vproof}{of Lemmas~\ref{lemC} and~\ref{lemmassene}.}
The proof of the embedding $W^{1,1}\big(I;X\big)\inj C_{\b,\vu}\big(\ovl I;X\big)$ is very standard and we omit its proof. Now, assume that $X\inj L^2(\Omega)$ with dense embedding. We infer that $L^2(\Omega)\inj X^\star.$ It follows that for any $v\in X,$
\begin{gather*}
\|v\|_{L^2(\Omega)}^2=(v,v)_{L^2(\Omega)}=\langle v,v\rangle_{L^2(\Omega),L^2(\Omega)}=\langle v,v\rangle_{X^\star,X}.
\end{gather*}
 We then note that $M\in C^1(\ovl I;\R),$ $E\in C^1(\ovl I;\R)$ and,
\begin{gather}
\label{demlemC1}
M(t)=M(s)+\int_s^t\big\langle u(\sigma),u^\p(\sigma)\big\rangle_{X,X^\star}\d\sigma,	\\
\label{demlemC2}
E(t)=E(s)+\int_s^t\big(-\Delta u(\sigma),u^\p(\sigma)\big)_{L^2(\Omega)}\d\sigma,
\end{gather}
for any $t,s\in\ovl I,$ as soon as $u\in\Dr(\ovl I;X),$ for \eqref{demlemC1} and $u\in\Dr(\ovl I;D(A)),$ for \eqref{demlemC2}. Applying Hölder's inequality in time and Young's inequality, one obtains,
\begin{gather}
\nonumber
\|u(t)\|_{L^2(\Omega)}^2\le\|u(s)\|_X\|u(s)\|_{X^\star}+\|u\|_{L^p(I;X)}^2+\|u^\p\|_{L^{p^\p}(I;X^\star)}^2,		\\
\label{demlemC3}
\|\nabla u(t)\|_{L^2(\Omega)}^2\le\|u(s)\|_{L^2(\Omega)}\|\Delta u(s)\|_{L^2(\Omega)}
					+\|\Delta u\|_{L^p(I;L^2)}^2+\|u^\p\|_{L^{p^\p}(I;L^2)}^2,
\end{gather}
for any $t,s \in\ovl I.$ Let $(I_n)_{n\in\N}\subset I$ be a increasing sequence (in the sense of the inclusion) of open bounded intervals such that $\bigcup_{n\in\N} I_n=I.$ Integrating in $s$ and applying, one more time, Hölder's and Young's inequalities, we have,
\begin{gather*}
|I_n|\,\|u\|_{C_\b(\ovl{I_n};L^2)}^2\le(1+|I_n|)\left(\|u\|_{L^p(I;X)}+\|u\|_{W^{1,p^\p}(I;X^\star)}\right)^2,
\end{gather*}
for any $n\in\N.$ Dividing by $|I_n|,$ letting $n\nearrow\infty$ and proceeding in the same way in \eqref{demlemC3}, we arrive at,
\begin{gather}
\label{demlemC4}
\|u\|_{C_\b(\ovl I;L^2)}\le(1+|I|^{-\frac12})\left(\|u\|_{L^p(I;X)}+\|u\|_{W^{1,p^\p}(I;X^\star)}\right),	\\
\label{demlemC5}
\|\nabla u\|_{C_\b(\ovl I;L^2)}\le(1+|I|^{-\frac12})\left(\|u\|_{L^p(I;D(A))}+\|u\|_{W^{1,p^\p}(I;L^2)}\right),
\end{gather}
with the convention $|I|^{-\frac12}=0,$ if $|I|=\infty.$ Since $X\inj X^\star$ and $D(A)\inj L^2(\Omega),$ we prove Lemma~\ref{lemC} by density with \eqref{demlemC4}--\eqref{demlemC5} (Lemma~\ref{lemden}). Finally, Lemma~\ref{lemmassene} is a consequence of \eqref{demlemC1}--\eqref{demlemC2} and Lemmas~\ref{lemden}--\ref{lemC}.
\medskip
\end{vproof}

\noindent
\textbf{Acknowledgements} \\
\baselineskip .5cm
The research of J.~I.~D{\'{\i}}az was partially supported by the project ref. MTM2017-85449-P of the DGISPI (Spain) and the Research Group MOMAT (Ref. 910480) of the UCM.

\baselineskip .4cm

\def\cprime{$'$}

\end{document}